\input amstex
\input epsf
\newdimen\unit
\unit=0.025in
\def\plot#1 #2 #3 {\rlap{\kern#1\unit\raise#2\unit\hbox{$#3$}}}
\magnification=\magstep1
\baselineskip=13pt
\documentstyle{amsppt}
\vsize=8.7truein
\CenteredTagsOnSplits
\NoRunningHeads
\def\today{\ifcase\month\or
  January\or February\or March\or April\or May\or June\or
  July\or August\or September\or October\or November\or December\fi
  \space\number\day, \number\year}
\def\dist{\operatorname{dist}}

\def\FF{\Cal{F}}
\def\Pr{\bold{P\ }}
\def\EE{\bold{E\ }}
\def\per{\operatorname{per}}
\NoRunningHeads
\topmatter
\title The Distance Approach to Approximate Combinatorial Counting \endtitle
\author Alexander Barvinok and  Alex Samorodnitsky \endauthor
\address Department of Mathematics, University of Michigan, Ann Arbor,
MI 48109-1109 \endaddress
\email barvinok$\@$math.lsa.umich.edu  \endemail
\address Institute for Advanced Study, 
Einstein Drive, Princeton, NJ 08540  \endaddress
\email asamor$\@$ias.edu \endemail
\date  May 2000 \enddate
\thanks The research of the first
author was partially supported by NSF Grant DMS 9734138.
\endthanks
\abstract We develop general methods to obtain 
fast (polynomial time) estimates of the cardinality of a combinatorially 
defined set via solving some randomly generated optimization problems
on the set.
Geometrically, we estimate the cardinality of a subset 
of the Boolean cube via the average distance from a point in the cube 
to the subset. As an application, we present a new randomized polynomial time
algorithm which approximates the permanent of a 0-1 matrix by solving 
a small number of Assignment problems. 
\endabstract
\endtopmatter
\document

\head 1. Introduction \endhead

A general problem of combinatorial counting can be stated as
follows: given a family $\FF \subset 2^X$ of subsets of the ground
set $X$, compute or estimate the cardinality $|\FF|$ of the family.
We would like to do the computation efficiently, {\it in polynomial 
time}. 
Of course,
one should clarify what ``given'' means, especially since in most 
interesting cases $|\FF|$ is exponentially large in the 
cardinality $|X|$ of the ground set. 
Following the earlier paper [Barvinok 97a], 
we assume that the family $\FF$ is defined by 
its {\it Optimization
Oracle}:
\specialhead (1.1) Optimization Oracle defining a family 
$\FF \subset 2^X$ \endspecialhead
\noindent{\bf Input:} A set of integer weights $\gamma_x: x \in X$. 
\medskip
\noindent{\bf Output:} The number 
$\displaystyle \min_{Y \in \FF} \sum_{x \in Y} 
\gamma_x$. 
\bigskip 
That is, for any given integer weighting 
$\{\gamma_x\}$ on the set $X$, we should be 
able to produce the minimum weight of a subset $Y \in \FF$. 
As is discussed in [Barvinok 97a], for many important families $\FF$ 
the Optimization Oracle is readily available.
The following example is central for this paper.   
\example{(1.2) Example: Perfect matchings in a graph}
Let $G=(V,E)$ be a graph with the set $V$ of vertices and 
set $E$ of edges. We assume that $G$ has no loops 
(edges whose endpoints coincide) and no isolated vertices.
 A set $M \subset E$ of edges is called a {\it matching}
in $G$ if every vertex of $G$ is incident to at most one edge from $M$.
A matching $M$ is called {\it perfect} if every vertex of $G$ is 
incident to precisely one edge from $M$.
Let $\FF \subset 2^E$ be the set of all 
{\it perfect matchings} in $G$.
The problem of computing or estimating $|\FF|$ efficiently 
is one of the hardest and most intriguing problems of combinatorial 
counting, see, for example, [Lov\'asz and Plummer 86], 
[Jerrum and Sinclair 89], [Jerrum 95] and [Jerrum and Sinclair 97].
 
We observe that Optimization Oracle 1.1 can be efficiently 
constructed. Indeed, if we assign integer weights $\gamma_e$: $e \in E$ 
to the edges of the graph, the minimum weight of a perfect matching 
can be computed in $O(|V|^3)$ time, see, for example Section 11.3 of
[Papadimitriou and Steiglitz 98]. 

A particularly interesting case is that of a {\it bipartite} graph $G$
when the vertices of $G$ are partitioned into two classes,
$V=V^+ \cup V^-$ such that every edge $e \in E$ has one endpoint in $V^+$ 
and the other in $V^-$. Then the number of perfect matchings 
in $G$ is equal to the permanent of a 0-1 matrix associated with $G$
(see also Section 5).
In this case, the corresponding optimization problem is known as 
the Assignment Problem. It is not 
only ``theoretically easy'', but in practice large instances are 
routinely solved as the Assignment Problem is a particular case of the minimum 
cost network flow problem (see, for example, Section 11.2 of 
[Papadimitriou and Steiglitz 98]).  
\endexample
Other interesting and generally difficult problems of 
combinatorial counting where the Optimization Oracle is provided 
by classical combinatorial optimization algorithms include counting 
bases in matroids, counting independent sets in 
matroids and counting bases in the intersection of two matroids 
over the same ground set, see [Jerrum and Sinclair 97] for 
a discussion of the counting problems and [Papadimitriou and Steiglitz 98]
for a description of the underlying optimization algorithms.
Particularly interesting special cases of those 
problems include counting spanning trees, counting forests and 
counting spanning subgraphs in a given graph and counting 
non-degenerate maximal minors in a given rectangular matrix over 
$GF(2)$.
Some of the problems, such as counting spanning trees, admit a simple 
and efficient solution, 
others, such as counting matchings of all sizes in a graph, are known 
to be hard to solve exactly but can be solved approximately and still others, 
such as counting bases in 
matroids, are solved only in special cases.
The problem of counting perfect matchings in a given graph, arguably 
the most famous problem of them all, still resists all attempts to solve it
in full generality (see also Section 5).

The most general approach to combinatorial counting has been via 
Monte Carlo method. The key component of the method is the ability 
to sample a random point from the (almost) uniform distribution on 
$\FF$. Often, to achieve this, 
a Markov chain on the set $\FF$ is generated, so 
that it converges rapidly to the uniform distribution on $\FF$ 
(see [Jerrum and Sinclair 97] for a survey). 
This approach resulted, for example,
in finding a polynomial time randomized algorithm to count matchings 
of all sizes in a given graph with a 
prescribed accuracy [Jerrum and Sinclair 89]. When the 
Markov chain approach works, it produces incomparably better 
results than the method of this paper. However, for many important 
counting problems, some of which are mentioned above, it is either not 
clear how to generate a rapidly mixing Markov chain or, when there is 
a ``natural'' candidate, it seems to be 
extremely hard to prove that the chain is indeed converging rapidly enough
to the steady state (cf. [Jerrum and Sinclair 97]). 
In contrast, our approach produces very crude bounds, but it is
totally insensitive to the fine structure of $\FF$, so it is ready 
to handle a broad class of problems. 
In [Barvinok 97a], it was shown that the method allows one 
to decide whether the size $|\FF|$ is exponentially large in the size $|X|$ 
of the ground set in some precisely defined sense. In this paper,
we improve the estimates of [Barvinok 97a] in several directions and 
apply them to new
problems, notably to the problem of estimating the permanent of 
a given 0-1 matrix.

The main idea of our approach is as follows. Given a family $\FF$, we 
identify it with a subset $F$ of a metric space $(\Omega, d)$, such that 
for any given point $x \in \Omega$ the distance 
$\displaystyle d(x, F)=\min_{y \in F} d(x,y)$ can be quickly computed 
using Optimization Oracle 1.1 for $\FF$. Then we estimate the 
cardinality $|F|$ from the distance $d(x, F)$ for a typical 
$x \in \Omega$. Intuitively, 
if $|F|$ is small, we expect the distance $d(x, F)$ from a random point 
$x \in \Omega$ to be large and vice versa. In this paper, 
$\Omega$ is the Boolean cube $\{0,1\}^n$ and $d$ is either the
Hamming distance or its modification, although as we discuss in 
Section 7, some other possibilities may be of interest.
Thus our approach can be considered as a refinement of the classical 
Monte-Carlo method: we do not only register how often a randomly sampled point 
$x \in \Omega$ lands in the target set $F$, but also take into account the  
distance $d(x,F)$. This allows us to get non-trivial bounds even 
when $|F|$ is exponentially small with respect to $|\Omega|$ so that $x$ 
typically misses $F$.  

The paper is organized as follows.

In Section 2, we introduce a  ``geometric cousin'' of Optimization
Oracle 1.1. Distance Oracle 2.2 describes a subset $F$ of the 
Boolean cube $\{0,1\}^n$ by computing a suitably 
defined distance $d$ from a given point in the cube to the set. We show
how to construct embeddings $\phi: \FF \longrightarrow \{0,1\}^n$, so that 
the Distance Oracle for the image $F=\phi(\FF)$ is derived from 
the Optimization Oracle for $\FF$. We show that in some important cases
(for example, when $\FF$ is the set of perfect matchings in a graph),
we can ``squeeze'' $\FF$ into a substantially smaller cube 
than we would have expected for a general family $\FF$.   

In Section 3, we describe the bounds obtained by choosing $d$ to be 
the Hamming distance in the cube. 
The bounds are sharp, meaning that we can't possibly 
estimate (in polynomial time) the cardinality of a subset 
$F \subset \{0,1\}^n$ better if the only information available is the 
Hamming distance from any given point $a \in \{0,1\}^n$ to the set $F$.
Remarkably, the lower and the upper bound for $\alpha=n^{-1} \log_2|F|$ 
converge when 
$\alpha \approx 0$ or $\alpha \approx 1$ and diverge the greatest 
when $\alpha=1/2$. 

In Section 4, we describe how to get better bounds for small sets by 
using a suitably defined ``randomized Hamming distance'', which
ignores a (random) part of the information contained in the standard 
Hamming distance. 
The isoperimetric problems arising here seem to be interesting in their 
own right. The 
proofs are not complicated but somewhat lengthy and therefore postponed till 
Section 6.

In Section 5, we apply our methods to get a new polynomial time 
algorithm to approximate the permanent of 
a given 0-1 matrix. 
Geometrically, we represent the set of the perfect 
matchings of the underlying bipartite graph on $n+n$ vertices as a
subset $F$ of the Boolean cube $\{0,1\}^m$ with $m=O(n \ln n)$ and estimate
$|F|$ from the Hamming distance of a random point in the cube to $F$. 
We find the distance in question by averaging solutions of some randomly 
generated Assignment problems. We compare our method with other algorithms 
available in the literature. In particular, we show that our method 
allows us to recognize $n \times n$ matrices whose permanents
are subexponential in $n$ (Corollary 5.4).

In Section 6, we supply proofs of the results of Section 4.  

In Section 7, we discuss possible ramifications of our approach and 
its relations with the Monte-Carlo method.

\head 2. Distance Oracle and Cubical Embeddings \endhead

The idea of our method is to represent $\FF$ geometrically as a 
subset $F$ of the Boolean cube and then derive estimates of $|\FF|$ 
using the average distance from a point in the cube to $F$. 
\definition{(2.1) Definitions} Let $C_n=\{0,1\}^n$ be the Boolean 
cube and let $\dist$ be the Hamming distance in $C_n$, that is 
$$\dist(a,b)=\sum_{i: \alpha_i \ne \beta_i} 1 \qquad 
\text{for} \quad a=(\alpha_1, \ldots, \alpha_n), 
\ b=(\beta_1, \ldots, \beta_n) 
\in C_n.$$
More generally, let us fix $n$ functions $d_i: \{0,1\} \times \{0,1\} 
\longrightarrow {\Bbb Z}$, $i=1, \ldots, n$, which we interpret as 
{\it penalties}. We assume that $d_i \geq 0$ and that 
$d(0,0)=d(1,1)=0$. Finally, 
let 
$$d(a,b)=\sum_{i=1}^n d_i(\alpha_i, \beta_i), \quad \text{where} \quad 
a=(\alpha_1, \ldots, \alpha_n) \quad \text{and} \quad 
b=(\beta_1, \ldots, \beta_n)$$
be the $L^1$ distance function determined by the penalties $\{d_i\}$.

If $d_i(\alpha, \beta)=1$ whenever $\alpha \ne \beta$ then 
$d(a,b)=\dist(a,b)$.
 
For a subset $B \subset C_n$ and a point $a \in C_n$, let 
$$d(a,B)=\min_{b \in B} d(a,b)$$
be the distance from $a$ to $B$. In particular,
let 
$$\dist(a,B)=\min_{b \in B} \dist(a,b)$$ 
be the Hamming distance from a point $a$ to the subset $B$.
\enddefinition 
We will be working with the following ``geometric cousin'' of 
Optimization Oracle 1.1.
\specialhead (2.2) Distance Oracle defining a set $F \subset C_n$ 
\endspecialhead
\noindent{\bf Input:} A point $a \in C_n$ and penalties $d_i: \{0,1\} 
\times \{0,1\} \longrightarrow {\Bbb Z}$, $i=1, \ldots, n$.
\medskip
\noindent{\bf Output:} The number $d(a,F)$.
\bigskip
There is an obvious way to associate with a family $\FF \subset 2^X$ a
subset $F \subset C_{|X|}$ of the Boolean cube.  
\subhead(2.3) Straightforward embedding \endsubhead 
Let us identify the ground set $X$ with the set $\{1, \ldots, n\}$,
$n=|X|$. Let $\FF$ be a family of subsets of $\{1, \ldots, n\}$ given by its 
Optimization Oracle.
For a subset $Y \in \FF$ let us
define the indicator $y \in C_n$, $y=(\eta_1, \ldots, \eta_n)$ 
by 
$$\eta_i=\cases 1 &\text{\ if\ } i \in Y \\ 0 &\text{\ if \ } i \notin Y.
 \endcases$$
Let $F=\bigl\{y \in C_n: Y \in \FF \bigr\}$ be the set of all indicators
of subsets $Y \in \FF$. 

Let us construct the Distance Oracle for the set $F \subset C_n$.
 Given a point $a=(\alpha_1, \ldots, \alpha_n) \in C_n$ and penalties
$d_i$, $i=1, \ldots, n$, let 
us define weights $\gamma_i$ by $\gamma_i=d_i(\alpha_i, 1)-d_i(\alpha_i,0)$.
Then for a set $Y \subset \{1, \ldots, n\}$ and its indicator
$y=(\eta_1, \ldots, \eta_n) \in C_n$, we have
$$\sum_{i \in Y} \gamma_i =\sum_{i \in Y} 
\bigl(d_i(\alpha_i, 1)-d_i(\alpha_i, 0)\bigr)=\sum_{i=1}^n d_i(\alpha_i, 
\eta_i)
-\sum_{i=1}^n d_i(\alpha_i, 0)=d(a,y)-d(a,0).$$
Hence, given the output 
$$\lambda=\min_{Y \in \FF} \sum_{i \in Y} \gamma_i$$ 
of Oracle 1.1 for the family $\FF$, we can easily compute the output 
$$d(a,F)=\lambda+d(a,0)$$ of Oracle 2.2 for the 
set $F$.    
Thus, given an Optimization Oracle 1.1 for a family $\FF \subset 2^X$,
we can efficiently construct a Distance Oracle 2.2 for 
a set $F \subset C_n$, $n=|X|$, such that $|F|=|\FF|$. 
\medskip
To be able to estimate the cardinality $|\FF|$ with a better precision,
we would like to embed $\FF$ into a smaller Boolean cube. Sometimes
this is indeed possible.
\subhead (2.4) Economical embedding \endsubhead Suppose that the
ground set $X$ can be represented as a union $X=X_1 \cup \ldots \cup X_k$ 
of (not necessarily disjoint) parts $X_i$, so that 
$|Y \cap X_i|=1$ for every subset $Y \in \FF$ and every $X_i$. In other
words, every member of $\FF$ is a transversal of the cover of $X$ 
by $X_1, \ldots, X_k$. 
Let 
$$m_i =\lceil \log_2 |X_i| \rceil \quad \text{and} \quad m=\sum_{i=1}^k m_i.$$
We construct an embedding $\FF \longrightarrow C_m$ as follows.

First, we index the elements of $X_i$ by distinct binary strings of 
length $m_i$,
that is, we choose an embedding $ \phi_i: X_i \longrightarrow C_{m_i}$.
Thus for any $x \in X_i$ the point $\phi_i(x)$ is a binary string of 
length $m_i$ and $\phi_i(x) \ne \phi_i(y)$ provided $x \ne y$.

Let us identify 
$$C_m= C_{m_1} \times \ldots \times C_{m_k}.$$
For a subset $Y \in \FF$, let us define $y \in C_m$ as
$$y=(y_1, \ldots, y_k), \quad \text{where} \quad 
y_i=\phi_i(Y \cap X_i) \in C_{m_i}.$$
Note that $y$ is well-defined, since every intersection $Y \cap X_i$ consists 
of a single point. 
Let $F=\bigl\{y \in C_m: Y \in \FF \bigr\}$. 
Clearly, $|F|=|\FF|$.

Given an Optimization Oracle 1.1 for $\FF$, let us construct a
Distance Oracle 2.2 for $F$. The input of Oracle 2.2 consists of 
a point $a \in C_m$ (binary string of length $m$) 
and penalty functions $\{d_i: i=1, \ldots m\}$.
We view
$a$ as  
$$a=(a_1, \ldots, a_k), \quad \text{where} \quad a_i \in C_{m_i}.$$
The penalties $d_i$, $i=1, \ldots, m$
 give rise to the $L^1$ distance function $d$ on binary 
strings, cf. Definition 2.1.
For a point $x \in X$, let us define its weight $\gamma_x$ by 
$$\gamma_x=\sum_{i: \ x \in X_i} d\bigl(a_i, \phi_i(x) \bigr). \tag2.4.1$$
Let $Y \in \FF$ be a set and let $y \in C_m$ be the point representing $Y$.
We observe that
$$\sum_{x \in Y} \gamma_x = \sum_{x \in Y} 
\sum_{i: \ x \in X_i} d\bigl(a_i, \phi_i(x) \bigr)=
\sum_{i=1}^k d(a_i, y_i)=d(a,y).$$
Hence, the outputs of Oracles 1.1 and 2.2 coincide: 
$$\min_{Y \in \FF} \sum_{x \in Y} \gamma_x=\min_{y \in Y} d(a,y).$$
Thus, given an Optimization Oracle 1.1 for a family $\FF \subset 2^X$,
we can efficiently construct a Distance Oracle 2.2 for a set $F \subset 
C_m$, such that $|F|=|\FF|$. More precisely, given a point 
$a \in C_m$ and penalties $\{d_i\}$, 
we compute weights $\{\gamma_x\}$ on $X$ by (2.4.1) in
$O(k|X| \ln |X|)$ time and then apply Optimization Oracle 1.1 to find 
the minimum weight $\lambda$ of a subset $Y \in \FF$ in this weighting. 
The distance $d(a, F)$ is equal to $\lambda$.
\example{(2.5) Example: Embedding perfect matchings}
Let $\FF$ be the family of all perfect matchings in a graph 
$G=(V,E)$, see Example 1.2. The straightforward embedding (2.3) identifies 
$\FF$ with a subset $F$ of the Boolean cube $\{0,1\}^{|E|}$ and provides 
us with Distance Oracle 2.2 for $F$. We will be better off using
the economical embedding (2.4). Indeed, for a vertex $v \in V$ of $G$, let 
$E_v$ be the set of edges of $G$ incident to $v$. Then 
$\displaystyle E=\cup_{v \in V} E_v$ and every perfect matching
has exactly one edge in every set $E_v$.
Hence the embedding (2.4) identifies $\FF$ with a subset $F$ of 
the Boolean cube $\{0,1\}^m$, where 
$$m=\sum_{v \in V} \lceil \log_2 |E_v| \rceil$$ and 
provides us with Distance Oracle 2.2 for $F$. 
Given a point $a \in C_m$, by (2.4.1) we compute weights $\gamma_e$ 
on the edges $E$ in $O(|E| \ln |E|)$ time (since every edge $e \in E$ 
belongs to exactly two sets $E_v$) and then find 
the minimum weight $\lambda$ of a perfect matching in $G$ in $O(|V|^3)$ time.
The distance $d(a, F)$ from $a$ to $F$ is equal to $\lambda$.  

Typically, if the graph has $|V|=n$ vertices and $\Omega(n^2)$ edges,
the dimension of the straightforward embedding will be $\Omega(n^2)$,
whereas the dimension of the economical embedding will be  
$O(n \ln n)$. We observe that for bipartite graphs we can reduce the 
dimension further by a factor of 2 at least by choosing
$$m=\min \Bigl\{\sum_{v \in V^+} \lceil \log_2 |E_v| \rceil, 
\sum_{v \in V^-} \lceil \log_2 |E_v| \rceil \Bigr\},$$
since every perfect matching $M \subset E$ will be a transversal 
of either partition $\displaystyle E=\cup_{v \in V^+} E_v$ or
$\displaystyle E=\cup_{v \in V^-} E_v$. 
\endexample   
Another natural case of economical embedding 2.4 arises when $\FF$ is the 
set of common bases of two matroids on the same ground set, one of which 
is a transversal matroid. It would be interesting to find out if 
similar economical embeddings can be constructed for a broader class 
of families $\FF \subset 2^X$, 
for example, when $\FF$ consists of ``small'' sets, 
that is, when $|Y| << |X|$ for any $Y \in \FF$.

\head 3. Estimating Cardinality from the Hamming Distance \endhead

In this section, we obtain estimates of the cardinality of 
a subset $F \subset C_n$ if we choose $d_i(0,1)=d_i(1,0)=1$,
$i=1, \ldots, n$ in Distance Oracle 2.2. In other words, 
we estimate $|F|$, provided we can compute the Hamming distance 
$\dist(x,F)$ to $F$ from any given point $x \in C_n$, cf. Definitions 2.1.
Our main tool is the {\it average} Hamming distance from a point to the 
set.     
\definition{(3.1) Definition} Let $A \subset C_n$ be a subset 
of the Boolean cube. Let 
$$\Delta(A)={1 \over 2^n} \sum_{x \in C_n} \dist(x, A)$$
be the average Hamming distance from a point in the cube to the set $A$.

Obviously, $\Delta(A) \leq \Delta(B)$ if $B \subset A$.
\enddefinition
\example{(3.2) Example: Set consisting of a single point}
Suppose that the set $A$ is a point. Without loss of generality 
we assume that $A=\{(0,\ldots, 0)\}$.
Then, for $x=(\xi_1, \ldots, \xi_n)$ we have $\dist(x,A)=\dist(x,0)=
\xi_1 + \ldots + \xi_n$ and
$$\Delta(A)={1 \over 2^n} \sum_{x \in C_n} \dist(x, A)=
{1 \over 2^n} \sum_{x \in C_n} (\xi_1+ \ldots + \xi_n)={n \over 2}.$$
It follows then that $\Delta(A) \leq n/2$ for any non-empty 
$A \subset C_n$ and that $\Delta(A)=n/2$ if and only if $A$ consists
of a single point. 
\endexample
Our first objective is to present a probabilistic algorithm that computes 
$\Delta(A)$ 
approximately by averaging $\dist(x, A)$ for a number of randomly 
chosen $x \in C_n$. 
\specialhead (3.3) Algorithm for computing $\Delta(A)$ \endspecialhead
\noindent{\bf Input:} A set $A \subset C_n$ defined by its 
Distance Oracle 2.2 and a number $\epsilon>0$.
\medskip
\noindent{\bf Output:} A number $\alpha$ approximating $\Delta(A)$ within 
error $\epsilon$.
\medskip
\noindent{\bf Algorithm:} Let $k=\lceil 48n/\epsilon^2 \rceil$. 
Sample $k$ points $x_1, 
\ldots, x_k \in C_n$ independently at random from the uniform 
distribution in the cube $C_n$. Apply Distance Oracle 2.2 to find
$\dist(x_i, A)$, $i=1, \ldots, k$. Compute 
$\displaystyle \alpha={1 \over k} \sum_{i=1}^k \dist(x_i, A)$. 
Output $\alpha$.
\bigskip
To prove that Algorithm 3.3 indeed approximates $\Delta(A)$ with the 
desired accuracy, we need a couple of technical results. The first lemma 
supplies us with important {\it concentration inequalities} for 
the Boolean cube.
\proclaim{(3.4) Lemma} Let $C_N=\{0,1\}^N$ be the Boolean cube and 
let $f: C_N \longrightarrow {\Bbb R}$ be a function such that 
$$|f(x)-f(y)| \leq \dist(x,y) \quad \text{for all} \quad x,y \in C_N.$$
Let 
$$\EE(f)={1 \over 2^N} \sum_{x \in C_N} f(x)$$ be 
the average value of $f$.
Let $\Pr$ denote the uniform probability measure on $C_N$, so 
$\Pr(A)=|A|/2^N$ for a set $A \subset C_N$.

Then for any $\delta>0$ 
$$\Pr\Bigl\{x \in C_N: |f(x)-\EE(f)| \geq \delta \Bigr\} \leq 
2 \exp\Bigl\{{-\delta^2 \over 16 N} \Bigr\}.$$
\endproclaim
\demo{Proof} See Sections 6.2 and 7.9 of [Milman and Schechtman 86].
{\hfill \hfill \hfill} \qed
\enddemo
The next lemma provides a useful ``scaling'' trick.
\proclaim{(3.5) Lemma} Let us fix positive integers $k$ and $n$ and let 
$N=kn$. Let us identify $C_N=C_n \times \ldots \times C_n=(C_n)^k$.
Thus a point $x \in C_N$ is identified with a $k$-tuple 
$x=(x_1, \ldots, x_k)$, where $x_i \in C_n$ for $i=1, \ldots, k$.

For a subset $A \subset C_n$, let 
$B=A\times \ldots \times A=A^k \subset C_N$. Then 
$$\dist(x, B)=\sum_{i=1}^k \dist(x_i, A) \quad \text{for any} \quad 
x=(x_1, \ldots, x_k) \in C_N$$ and 
$$\Delta(B)=k\Delta(A).$$
\endproclaim
\demo{Proof} 
Clearly,
$$\dist(x,y)=\sum_{i=1}^k \dist(x_i, y_i) \quad \text{for all} \quad
x, y \in C_N,$$
hence the first identity follows.
Next,
$$\split \Delta(B)&={1 \over 2^N} \sum_{x \in C_N} \dist(x,B)=
{1 \over 2^N} \sum_{x_1, \ldots, x_k \in C_n} 
\sum_{i=1}^k \dist(x_i, A)\\ 
&={k2^{n(k-1)} \over 2^{nk}} \sum_{x \in C_n} 
\dist(x, A)={k \over 2^n} \sum_{x \in C_n} \dist(x,A)=k\Delta(A).
\endsplit$$
{\hfill \hfill \hfill} \qed
\enddemo
Now we can prove correctness of Algorithm 3.3. 
\proclaim{(3.6) Theorem} With probability at least 0.9, the output 
$\alpha$ of Algorithm 3.3 satisfies the inequality 
$|\Delta(A)-\alpha| \leq \epsilon$.
\endproclaim
\demo{Proof} Let $N=nk$ and let us identify
$C_N=(C_n)^k$ as in Lemma 3.5. Let $B=A^k \subset C_N$.
Let $f: C_N \longrightarrow {\Bbb R}$ be defined by 
$f(x)=\dist(x,B)$. Applying Lemma 3.4 with $\delta=k\epsilon$ and 
observing that $\EE(f)=\Delta(B)$,
we conclude that  
$$\Pr\Bigl\{x: |\dist(x,B)-\Delta(B)| \geq k\epsilon \Bigr\} \leq 
2 \exp\Bigl\{-{(\epsilon k)^2 \over 16N} \Bigr\}=
2 \exp\Bigl\{-{\epsilon^2 k \over 16n} \Bigr\} \leq 0.1.$$
Since by Lemma 3.5
$$\Delta(B)=k\Delta(A) \quad \text{and} \quad 
{1 \over k} \sum_{i=1}^k \dist(x_i, A)={1 \over k} \dist(x, B)$$ 
for $x=(x_1, \ldots, x_k)$, we conclude that
$$\split 
&\Pr\Bigl\{x_1, \ldots, x_k: \Big|{1 \over k} \sum_{i=1}^k \dist(x_i, A)-
\Delta(A) \Big| \geq \epsilon \Bigr\}= \\ 
&\Pr\Bigl\{x: |\dist(x,B)-\Delta(B)| \geq k\epsilon \Bigr\} \leq 0.1,
\endsplit$$
and the proof follows.
{\hfill \hfill \hfill} \qed
\enddemo
\remark{Remark} Hence to evaluate $\Delta(A)$ 
within error $\epsilon$ we have to average $O(n \epsilon^{-2})$ values 
$\dist(x_i, A)$. By doing that, we allow probability 0.1 of failure. 
As usual, to attain a lower probability $\delta>0$ of 
failure, one should run Algorithm 3.3 $O(\ln \delta^{-1})$ times and then 
select the median of the computed $\alpha$'s (cf. [Jerrum {\it et al.} 86]).
For all applications, choosing $\epsilon=1$ will suffice and in many
cases $\epsilon=\sqrt{n}$ will do (cf. Section 5 and [Barvinok 97a]). 
Hence, often we 
will have to apply Oracle 2.2 only a constant number of times.  
\endremark
We would like to relate the value of $\Delta(A)$ to the cardinality 
$|A|$. 
\definition{(3.7) Definition. Entropy Function} For $0\leq x \leq 1/2$ let 
$$H(x)=x\log_2 {1 \over x} + (1-x) \log_2 {1 \over 1-x}.$$ 
We agree that $H(0)=0$.
Thus $H$ is an increasing concave function on the interval $[0,1/2]$.

We use the following estimate 
(see, for example, Theorem 1.4.5 of [van Lint 99])
$$\sum_{k=0}^r {n \choose k} \leq 2^{n H(r/n)} \qquad \text{for} \quad 
r \leq n/2. \leqno(3.7.1)$$
Also, we remark that around $x=+0$ we have 
$$H(x)=x \log_2 {1 \over x} +O(x) \quad \text{and} \quad 
H\Bigl({1 \over 2} -x \Bigr)=1-{2 \over \ln 2} x^2 +O(x^3) \leqno(3.7.2)$$
\enddefinition 
We will use the classical isoperimetric inequality for the Boolean cube
(see, for example, [Leader 91]).
\proclaim{(3.8) Harper's Theorem} Let $A \subset C_n$ be a set such 
that $$|A| \geq \sum_{k=0}^r {n \choose k}$$
for some integer $r$. Then, for any non-negative 
integer $t$ 
$$\big|\bigl\{x \in C_n: \dist(x,A) \leq t \bigr\} \big| \geq \sum_{k=0}^{r+t}
{n \choose k}.$$
\endproclaim
We are going to obtain an estimate of the cardinality of a set 
$A \subset C_n$ in terms of the average Hamming distance $\Delta(A)$ from
a point $x \in C_n$ to $A$. It is convenient to express the estimate 
in terms of a related quantity 
$$\rho=\rho(A)={1 \over 2} -{\Delta(A) \over n}.$$
As follows from Example 3.2, for every non-empty set $A \subset C_n$ we 
have $0 \leq \rho(A) \leq 1/2$. We observe that $\rho(A)=0$ if and only if 
$A$ consists of a single point and that $\rho(A)=1/2$ if and only if
$A$ is the whole cube $C_n$.   
\proclaim{(3.9) Theorem} Let $A \subset C_n$ be a non-empty set.
Let $$\rho={1 \over 2} - {\Delta(A) \over n}.$$ Then 
$$1-H\Bigl({1 \over 2}-\rho\Bigr) \leq {\log_2 |A| \over n} \leq H(\rho).$$
\endproclaim
Before we proceed with a formal proof, we would like to highlight 
some ideas.
\subhead (3.10) The idea of the proof. Extremal sets \endsubhead
Let $A \subset C_n$ be a set.
Concentration inequalities (Lemma 3.4) imply that the average 
distance $\Delta(A)$ is approximately equal to the distance
$\dist(x, A)$ from a ``typical'' point $x \in C_n$ to $A$.
For a given positive integer $t$, let us consider 
the $t$-neighborhood $A_t=\bigl\{x \in C_n: \dist(x, A) \leq t \bigr\}$ 
of $A$. We expect that $\Delta(A) \approx t_1$, where $t_1$ is the 
smallest value of $t$ such that $A_t$ covers ``almost all'' cube 
$C_n$. The neighborhood $A_t$ grows the slowest when 
$A$ is a ball in the Hamming metric, that is 
when $A=\bigl\{x: \dist(x,x_0) \leq r \bigr\}$ for some $x_0 \in C_n$ and 
some $r>0$, as follows from Harper's Theorem 3.8, cf. also [Leader 91].
Hence the upper bound for $n^{-1} \log_2 |A|$  
in Theorem 3.9 is attained (up to an $O(n^{-1/2})$
error term) when $A$ is a ball.  
The neighborhood $A_t$ grows the fastest when the points of 
$A$ are spread around in $C_n$. In any case, the size $|A_t|$ does not 
exceed the sum of sizes of the balls of radius $t$ centered at the 
points of $A$. Thus the lower bound for $n^{-1} \log_2 |A|$ 
in Theorem 3.9 is obtained 
from this ``packing'' type argument. One can show that if the points 
of $A$ are chosen at random in $C_n$, then with high probability      
the lower bound is indeed attained asymptotically. More precisely, 
let us fix a number $0<\beta <1$ and let $A$ be the set of 
$\lfloor 2^{\beta n} \rfloor$ points chosen at random from $C_n$.
Then with the probability that tends to 1 as $n$ grows to infinity,
$\beta=1-H\bigl({1 \over 2}-\rho\bigr)+O(n^{-1/2})$. 
The proof is straightforward, but 
technical and therefore omitted. 

Finally, we note that using average distance $\Delta(A)$ and the 
scaling trick (Lemma 3.5) allows us to get rid of  
$O(n^{-1/2})$ error terms in the proof.
\demo{Proof of Theorem 3.9}
 Let us choose a positive even integer $m$, let $N=mn$ and let 
us identify $C_N=(C_n)^m$, as in Lemma 3.5. 
Let $B=A^m \subset C_N$. Let us fix the uniform probability measure 
$\Pr$ on $C_N$.

Let $\alpha=\log_2|A|/n$, so $|A|=2^{\alpha n}$ and $|B|=2^{\alpha N}$.
Let $0 \leq \gamma \leq 1/2$ be a number such that $H(\gamma)=\alpha$ and 
let $r=\lfloor N \gamma \rfloor$. Then by (3.7.1) 
$$|B|=2^{N \cdot H(\gamma)} \geq \sum_{k=0}^r {N \choose k}.$$ 
Then Theorem 3.8 implies that  
$$\big| \bigl\{x \in C_N: \dist(x, B) \leq N/2-r \bigr\} 
\big| \geq \sum_{k=0}^{N/2}
{N \choose k}=2^{N-1}.$$
Therefore,
$$\Pr\Bigl\{x \in C_N: \dist(x, B) \leq {N \over 2}-r \Bigr\} 
\geq {1 \over 2}.$$
We have that $x=(x_1, \ldots, x_m)$ for some $x_i \in C_n$ and 
that $\dist(x,B)=\dist(x_1, A) + \ldots + \dist(x_m, A)$ (see Lemma 3.5).
Therefore, 
$$\Pr\Bigl\{(x_1, \ldots, x_m): {1 \over m} \sum_{i=1}^m \dist(x_i, A) 
\leq {N \over 2m} -{r \over m} \Bigr\} \geq {1 \over 2}. \tag1$$
Now we observe that 
$${N \over 2m} -{r \over m} \longrightarrow {n \over 2}-n \gamma 
\quad \text{as}  \quad 
m \longrightarrow +\infty. \tag2$$
Furthermore, by the Law of Large Numbers,
$${1 \over m}\sum_{i=1}^m \dist(x_i, A)
 \longrightarrow \Delta(A) \quad
\text{in probability} \qquad \text{as} \quad m \longrightarrow +\infty.
\tag3$$
Hence the assumption that $\Delta(A)> n/2 - n\gamma$ would contradict
(1)--(3). Thus we must have $\Delta(A) \leq n/2 -n\gamma$, which 
implies that $\gamma \leq \rho(A)$.
Hence $\alpha=H(\gamma) \leq H(\rho)$ and the upper bound is proven.

Let us prove the lower bound. We observe that for every point 
$b \in C_N$ and any $N/2 \geq s \geq 0$
$$\big|\{x \in C_N: \dist(x,b) \leq s \}\big| =\sum_{k=0}^s {N \choose k} \leq 
2^{N \cdot H(s/N)}.$$
Therefore,
$$\big|\{x \in C_N: \dist(x, B) \leq s \}\big|  \leq |B|
2^{N \cdot H(s/N)}=2^{N \cdot (H(s/N)+\alpha)}.$$
Hence
$$\Pr\bigl\{x \in C_N: \dist(x,B) \leq s\bigr\} 
\leq 2^{N \cdot (H(s/N)+\alpha -1)}.$$  
Therefore,
$$\Pr\Bigl\{(x_1, \ldots, x_m):
{1 \over m} \sum_{i=1}^m \dist(x_i, A) \leq s/m \Bigr\} \leq 
2^{N \cdot (H(s/N)+\alpha -1)}. \tag4$$
If $\Delta(A)=n/2$ then $A$ is a point and the lower bound in Theorem 3.9
is satisfied.
Otherwise, let us fix an $\epsilon>0$
such that $(1+\epsilon) \Delta(A)/n < 1/2$
and let 

$s=\lceil m(1+\epsilon) \Delta(A) \rceil$.
We have 
$$s/m \longrightarrow (1+\epsilon) \Delta(A) \quad \text{and} \quad 
s/N \longrightarrow (1+\epsilon) \Delta(A)/n \quad \text{as} \quad 
m \longrightarrow +\infty. \tag5$$
Hence the assumption that 
$\displaystyle H\bigl((1+\epsilon)\Delta(A)/n \bigr)+\alpha-1 < 0$
would contradict (3)--(5). Therefore,  
$\displaystyle H\bigl((1+\epsilon)\Delta(A)/n \bigr)+\alpha-1 \geq 0$ 
for any $\epsilon>0$ and $H\bigl(\Delta(A)/n \bigr)+\alpha-1 \geq 0$.
Since $\Delta(A)/n=0.5-\rho$, the proof
follows.
{\hfill \hfill \hfill} \qed
\enddemo
For applications, the most interesting case is when 
$\displaystyle n^{-1} \log_2 |A|$ is small, that is $\rho \approx 0$.
\proclaim{(3.11) Corollary} There exist positive constants $c_1$ and $c_2$
such that for any non-empty set $A \subset C_n$ and for 
$\displaystyle \rho={1 \over 2} -{\Delta(A) \over n}$ we have
$$c_1 \cdot \rho^2 \leq {\ln |A| \over n} \leq c_2 \cdot \rho 
\ln {1 \over \rho}.$$
In particular, for any $c_1 < 2$ and any $c_2>1$,
the inequality holds in a sufficiently small neighborhood of $\rho=0$.
\endproclaim
\demo{Proof} Follows from Theorem 3.9 by (3.7.2).
{\hfill \hfill \hfill} \qed
\enddemo
\subhead (3.12) Discussion \endsubhead
Figure 1 depicts the feasible region for $n^{-1} \log_2 |A|$ as described 
by Theorem 3.9. 
$$\hbox to 2.0 in{
\plot -30 0 {\epsffile{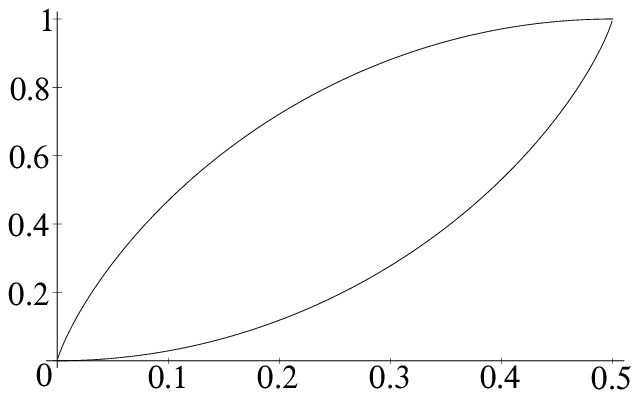}} 
\plot 100 30 {\rho} 
\plot 5 80 {\displaystyle {\log_2 |A| \over n}} 
\plot 35 75 {H(\rho)} 
\plot 60 40 {1-H\Bigl({1 \over 2}-\rho \Bigr)} 
\plot 35 10 {\bold{Figure \ 1}} 
\hfil}$$
Thus possible values of $n^{-1}\log_2 |A|$ 
with the given value of 
$\rho$ form a vertical interval between the two curves.
As we discussed in Section 3.10, asymptotically both bounds are sharp. 
Remarkably, the bounds converge at $\rho=0$ and $\rho=0.5$. On the 
other hand, the difference is the greatest when $\rho=1/4$. Thus, if 
the average Hamming distance from a point $x \in C_n$ to a set 
$A \subset C_n$ is $n/4$, the set $A$ can contain as many as 
$\displaystyle 2^{0.811n}$ points 
and as few as $\displaystyle 2^{0.189n}$ points.

Corollary 3.11 (with somewhat weaker constants and stated in different 
terms) 
together with the observation that the distance
$\dist(x,A)$ for a randomly chosen point $x \in C_n$ allows one to estimate
$\rho$ up to an $O(n^{-1/2})$ error constitute the main result 
of the earlier paper [Barvinok 97a]. 
Consequently, the main conclusion of [Barvinok 97a] is equivalent to 
stating that the Hamming distance 
to $A$ from a random point $x$ in the Boolean cube allows one to 
decide whether $|A|$ is exponentially large in $n$.
Theorems 3.6 and 3.9 make improvements of two kinds. First, we obtain 
sharp bounds valid for all $0 \leq \rho \leq 1/2$,
and second, by 
averaging several random distances (see Algorithm 3.3 and Theorem 3.6)
we get rid of the $O(n^{-1/2})$ error term. This allows us to obtain 
meaningful cardinality estimates for really small sets. For example,
if $A \subset C_n$ is a set such that 
$n^{-1} \log_2 |A| \sim n^{-\alpha}$, for some 
$0 < \alpha < 1$, by 
applying Algorithm 3.3 to approximate 
$\Delta(A)$ and Theorem 3.9 to interpret the results,
 the worst lower bound we can get for 
$n^{-1} \log_2 |A|$ is
$\sim n^{-2\alpha} \ln^{-2} n$ (this happens when $A$ is a ball in the 
Hamming metric, but we think it is a ``random set'', see 
Section 3.10) 
and the worst upper bound is   
$\displaystyle  \sim n^{-\alpha/2} \ln n$ (this happens when $A$ is 
a ``random set'' but we think that it is a ball).
Curiously, we can even distinguish in polynomial time between a set
consisting of a single point ($\rho=0$) and a set having more than one point
(one can show that $\rho \geq c/n$ for some $c>0$ in that case),
although apparently we can't distinguish between sets consisting of 
2 and 3 points respectively.

As we remarked earlier, in applications the value of 
$n^{-1} \log_2 |A|$ is usually small (cf. Examples 1.2 and 2.5).
Therefore, it 
is of interest to tighten the bounds for such sets. In the next 
section, we show that this is indeed possible: we demonstrate how to 
modify the definition of $\rho$, so that it remains efficiently computable 
and so that
$$c_3\cdot \rho^2 \ln {1 \over \rho} \leq {\ln |A| \over n} \leq 
c_4 \cdot \rho \ln {1 \over \rho}$$
for some $c_3, c_4>0$, which improves 
the inequality of Corollary 3.11 in the neighborhood of 
$\rho=0$.

\head 4. Randomized Hamming Distance \endhead 

Let us fix a number $0<p \leq 1$ and let $q=1-p$. In this section, 
we construct a quantity $\Delta(A,p)$, which measures the cardinality 
of ``small'' subsets $A \subset C_n$ of the Boolean cube 
in a somewhat more precise way than
the average Hamming Distance $\Delta(A)$ discussed in Section 3. In fact, 
$\Delta(A,1)=\Delta(A)$, so $\Delta(A)$ is a particular case of $\Delta(A,p)$.
\definition{(4.1) Definitions} 
Let $\Lambda_n$ be a copy of the Boolean cube $\{0,1\}^n$.  
We make $\Lambda_n$ a probability space by letting
$$\Pr\{l\}=p^{|l|} q^{n-|l|}, \quad \text{where} \quad 
|l|=\lambda_1 + \ldots + \lambda_n \quad \text{for} \quad 
l=(\lambda_1, \ldots, \lambda_n).$$  
Hence a vector $l=(\lambda_1, \ldots, \lambda_n)$ from $\Lambda_n$
is interpreted as a realization of $n$ independent random variables 
$\lambda_i$ such that 
$\Pr\{\lambda_i=1\}=p$ and $\Pr\{\lambda_i=0\}=q$. 

For $x, y \in C_n$ and an $l \in \Lambda_n$, where 
$x=(\xi_1, \ldots, \xi_n)$, $y=(\eta_1, \ldots, \eta_n)$ and 
$l=(\lambda_1, \ldots, \lambda_n)$, let 
$$d_l(x,y)=\sum_{i: \xi_i \ne \eta_i} \lambda_i.$$
In other words, we count disagreement in the $i$-th coordinate of 
$x$ and $y$ if and only if the value
of $\lambda_i$ is 1. 
Thus if $l=(1, \ldots, 1)$, we have $d_l(x,y)=\dist(x,y)$, the usual 
Hamming distance. 

For $l \in \Lambda_n$ and a set $A \subset C_n$, let
$$d_l(x,A)=\min_{y \in A} d_l(x,y).$$

Finally, let 
$$\Delta(A,p)= 
\sum_{l \in \Lambda_n} \sum_{x \in C_n}
d_l(x,A) {p^{|l|} q^{n-|l|} \over 2^n}.$$
In other words, $\Delta(A,p)$ is the expected value of $d_l(x,A)$, where 
$x=(\xi_1, \ldots, \xi_n)$ and $l=(\lambda_1, \ldots, \lambda_n)$ are vectors
of independent random variables such that $\Pr\{\lambda_i=1\}=p$, 
$\Pr\{\lambda_i=0\}=q$ and $\Pr\{\xi_i=0\}=\Pr\{\xi_i=1\}=1/2$.
Obviously, $\Delta(A,p) \leq \Delta(B,p)$ if $B \subset A$.

It follows that for a fixed non-empty $A \subset C_n$, the value 
$\Delta(A,p)$ is a polynomial in $p$ of degree at most $n$. 
\enddefinition
\example{(4.2) Example: Set consisting of a single point} 
Suppose that the set $A$ consists of a single point. Without loss
of generality we assume that $A=\{(0, \ldots, 0)\}$. 
Then for $x=(\xi_1, \ldots, \xi_n)$ 
and $l=(\lambda_1, \ldots, \lambda_n)$, 
$$d_l(x,A)=\sum_{i=1}^n \lambda_i \xi_i.$$
Interpreting $\lambda_i$ and $\xi_i$, $i=1, \ldots, n$ as independent 
random variables such that $\Pr\{\xi_i=1\}=\Pr\{\xi_i=0\}=1/2$ and 
$\Pr\{\lambda_i=1\}=p$, $\Pr\{\lambda_i=0\}=q$, we get 
$$\Delta(A,p)=\EE \sum_{i=1}^n \lambda_i \xi_i =\sum_{i=1}^n 
(\EE \lambda_i) (\EE \xi_i)={np \over 2}.$$ 
\endexample    
It follows then that for any non-empty set $A \subset C_n$ 
we have $\Delta(A,p) \leq np/2$ and that $\Delta(A,p)=np/2$ if and 
only if $A$ consists of a single point (we agreed that $p>0$).

As was the case with $\Delta(A)$, the functional $\Delta(A,p)$ can be 
easily computed by averaging. For a set $A \subset C_n$ 
defined by its Distance Oracle 2.2 and any $l=(\lambda_1, \ldots, \lambda_n)$ 
the value of $d_l(x,A)$ is computed by choosing the 
penalties $d_i(0,1)=d_i(1,0)=1$ when $\lambda_i=1$ and 
$d_i=0$ when $\lambda_i=0$.  

\specialhead{(4.3) Algorithm for Computing $\Delta(A,p)$}\endspecialhead 

\noindent {\bf Input:} A set $A \subset C_n$ given by its Distance Oracle 2.2,
a number $1 \geq p > 0$ and an $\epsilon>0$.
\medskip
\noindent{\bf Output:} A number $\alpha$ approximating $\Delta(A,p)$
within error $\epsilon$.
\medskip
\noindent{\bf Algorithm:} Let $k=\lceil 64n/\epsilon^2 \rceil$.
Sample $k$ points $x_1, \ldots, x_k \in C_n$ independently at random from 
the uniform distribution in $C_n$ and $k$ points 
$l_1, \ldots, l_k \in \Lambda_n$ 
independently at random from the distribution in $\Lambda_n$. Apply 
Distance Oracle 2.2 to compute $d_{l_i}(x_i, A)$, $i=1, \ldots, k$.
Compute $\displaystyle \alpha={1 \over k} \sum_{i=1}^k \dist_{l_i}(x_i, A)$.
Output $\alpha$.
\proclaim{(4.4) Theorem} With probability at least 0.9, the output 
$\alpha$ of Algorithm 4.3 satisfies the inequality 
$|\Delta(A, p)-\alpha| \leq \epsilon$.
\endproclaim
We postpone the proof till Section 6.

We are going to obtain estimates of the cardinality $|A|$ of a set 
$A \subset C_n$ in terms of the quantity $\Delta(A,p)$. As in 
Section 3, it is convenient to work with a related quantity
$$\rho=\rho(A,p)={p \over 2} - {\Delta(A,p) \over n}.$$
From Definitions 4.1, for any non-empty $A \subset C_n$, the function
$\rho(A,p)$ is a polynomial in $p$ of degree at most $n$. 
As follows from Example 4.2, $0 \leq \rho \leq p/2$ for any non-empty set 
$A \subset C_n$. Our estimate will be useful for ``small'' sets $A$ 
where $n^{-1} \ln |A|$ is close to 0.
\proclaim{(4.5) Theorem} Let $A \subset C_n$ be a non-empty set.
Let 
$$\rho={p \over 2} - {\Delta(A,p) \over n}.$$
Then
$${\rho^2 \over p} \leq {\ln |A| \over n}. \leqno(4.5.1)$$
Suppose that $\rho \leq 1/4$ and that 
$$p \geq {\ln 2 +\ln(1-2\rho) \over \ln(1-2\rho) -\ln(2\rho)}. \leqno(4.5.2)$$
Then 
$${\ln |A| \over n} \leq 2\rho \ln {1 \over 2 \rho} +(1-2\rho) 
\ln{1 \over 1-2\rho}.
\leqno(4.5.3)$$
\endproclaim
As we remarked earlier, the case interesting for applications is when 
$|A|$ is small, meaning that $n^{-1} \ln |A| \approx 0$.
\proclaim{(4.6) Corollary} Let us choose any
$c_3<1/(\ln 2) \approx 1.44$ and any $c_4>2$. Then
there exists a $\delta>0$ such that for any non-empty $A \subset C_n$ 
with $n^{-1} \ln |A| \leq \delta$ there exists a $0<p \leq 1$ such that 
for $\displaystyle \rho={p \over 2} -{\Delta(A,p) \over n}$ one has
$$c_3 \cdot \rho^2 \ln {1 \over \rho} \leq {\ln |A| \over n} 
\leq c_4 \cdot \rho \ln{1 \over \rho}.$$
\endproclaim
\demo{Proof}  
By (4.5.1), $\rho \leq \sqrt{n^{-1} \ln |A|} \leq \sqrt{\delta}$,
so $\rho(A,p)$ is small if $\delta$ is small, no matter what $p$ is.
We observe that for small positive $\rho$ the right hand side 
of (4.5.2) is of the order $(\ln 2) \ln^{-1} (1/\rho)$ and the right hand 
side of (4.5.3) is of the order $2 \rho \ln (1/\rho)$.

Given $c_3 < (\ln 2)^{-1}$ and $c_4>2$, 
let us choose $1/16>\delta>0$ in such a way that 
the right hand side of (4.5.2) does not exceed $(c_3)^{-1} \ln^{-1} (1/\rho)$ 
and the right hand side of (4.5.3) does not exceed $c_4 \rho \ln (1/\rho)$
for all $0< \rho < \sqrt{\delta}<1/4$.

We recall that $|A|=1$ if and only if $\rho=0$, in which case the 
bounds of Corollary 4.6 are satisfied by default.
Given a set $A \subset C_n$, $|A|>1$, 
let us choose the smallest $p \geq  0$ that 
satisfies the inequality (4.5.2). Then $0<p < 1$ since the 
right hand side of (4.5.2) is bounded below from 0 as 
a function of $p$ and $n$ and smaller than 1 for $0<\rho < 1/4$.  
Since $\rho(A,p)$ depends continuously on $p$, we must 
have equality in (4.5.2) (otherwise, we could have taken a smaller $p$). 
Thus $p \leq (c_3)^{-1} \ln^{-1} (1/\rho)$ and the proof 
follows by (4.5.1)--(4.5.3).
{\hfill \hfill \hfill} \qed
\enddemo 
\subhead (4.7) Extremal sets \endsubhead
Let us fix a $0<p \leq 1$ and an $\epsilon>0$.
Then there exists an $\alpha=\alpha(p, \epsilon)>0$ 
with the following property: if $A \subset C_n$ is 
a set of $\lfloor 2^{\alpha n} \rfloor$ points randomly chosen from 
the Boolean cube, then with the probability that tends to 1 as 
$n$ grows to infinity,
$n^{-1} \ln |A| < (2+\epsilon) \rho^2/p$.
Hence for any $p>0$ the bound (4.5.1) is tight up to a constant 
factor for sufficiently small random sets. The proof is rather technical
and therefore omitted. 

One can show that as long as $p$ 
satisfies (4.5.2), the bound (4.5.3) is asymptotically attained 
 on small {\it faces} of the cube $C_n$. 
Let us fix a $\delta>0$ (to be adjusted later), 
let $m=\lfloor \delta n \rfloor$ and 
let $A \subset C_n$ be an
$m$-dimensional face of the 
Boolean cube:
$$A=\Bigl\{(\xi_1, \ldots, \xi_n): \xi_i=0 \quad \text{for} \quad 
i=m+1, \ldots, n\Bigr\}.$$
Thus $|A|=2^m$. Moreover, a computation similar to that of
Example 4.2 shows that $\rho(A,p)=pm/2n$.
Hence we have 
$${\ln |A| \over n} = {2 \ln 2 \over p } \rho(A,p).$$ 
We observe that $\rho(A,p) \leq \delta/2$.
Hence for any small $\epsilon>0$ one can find 
$\delta=\delta(\epsilon)>0$ such that there exists $p$ 
satisfying (4.5.2) and such that  
$p <  (1+\epsilon) (\ln 2) \ln^{-1} (1/\rho)$. For such a $p$, we have 
$${\ln |A| \over n} \geq  {2 \over 1+\epsilon} \rho \ln{1 \over \rho},$$
so the bound (4.5.3) is indeed asymptotically tight for small sets.

Apparently, the sets $A$ having the largest cardinality among all sets 
with the given value of $\rho(A,p)$ evolve from the balls in the Hamming 
metric for $p=1$ (see Section 3.10) to faces at $p \longrightarrow 0$. 
Since faces are packed somewhat less tightly than balls,
we gain in Corollary 4.6 as compared to Corollary 3.11. 
\medskip
The proof of Theorem 4.5 is postponed till Section 6. 

Corollary 4.6 implies for small sets $A$ by ``tuning up'' $p$ we can 
get an additional logarithmic factor which brings the lower bound 
for $n^{-1} \ln |A|$ a little closer to the upper bound compared 
to the bound of Corollary 3.11. Any $p$ which is only slightly bigger 
than the bound (4.5.2) will do. 
For example,
if $A \subset C_n$ is a set such that $n^{-1} \ln |A|
\sim n^{-\alpha}$ for some $0<\alpha<1$, it follows by (4.5.1) that 
$\rho(A,p)=O(n^{-\alpha/2})$ for any $p$. Then we can choose
some $p=O(\ln^{-1} n)$ that satisfies (4.5.2) (
a particular suitable value of $p$ can be found, for 
example, by dichotomy).
Applying Algorithm 4.3
to approximate $\Delta(A,p)$ and Theorem 4.5 to interpret the 
results, for $n^{-1} \ln |A| $ we would obtain a lower 
bound of the form $\sim n^{-2 \alpha}/ \ln n$ at worst
and 
an upper bound of the form $\sim n^{-\alpha/2} \sqrt{\ln n}$ at worst, 
which is 
somewhat better than the bounds that could possibly 
be obtained by using the standard Hamming distance, see Section 3.12.

We are not going to use $\Delta(A,p)$ in what follows, but we find it 
interesting that some improvement in the cardinality estimate can be 
achieved by simply ignoring a (random) part of the information contained
in the standard Hamming distance.
   
\head 5. Application: Approximating the Permanent of a 0-1 Matrix \endhead

Let $A=(a_{ij})$ be an $n \times n$ matrix. The {\it permanent} of $A$ 
is defined by the expression 
$$\per A=\sum_{\sigma \in S_n} \prod_{i=1}^n a_{i \sigma(i)},$$
where $S_n$ is the symmetric group of all substitutions of the set 
$\{1, \ldots, n\}$. 
If $a_{ij} \in \{0,1\}$ for all $i$ and $j$ then $\per A$ counts 
perfect matchings in a bipartite graph $G_A=(V,E)$, constructed as 
follows. Let $V=V^+ \cup V^-$ be the set of vertices, where 
$V^+=\{1^+, \ldots, n^+\}$ and $V^-=\{1^-, \ldots, n^-\}$, and let 
$e=(i^+,j^-)$ be an edge of $G_A$ if and only if $a_{ij}=1$. Then 
$\per A$ is equal to the number of perfect matchings in $G_A$, cf. Example 
1.2. The problem of computing $\per A$ is $\#$ P-hard [Valiant 79] and
polynomial time algorithms for computing $\per A$ exactly are known 
only in few particular cases. For example, if the graph $G_A$ is 
planar (see [Lov\'asz and Plummer 86]), or more generally, has the genus 
bounded by some absolute constant [Gallucio and Loebl 99] then $\per A$ can be 
computed in polynomial time. If the permanent of a 0-1 
matrix is small (bounded by a polynomial in the size $n$ of the matrix), 
it can be computed in polynomial time, see Section 7.3 of [Minc 78] and 
[Grigoriev and Karpinski 87]. Finally,
the permanent of matrices (real or complex) of a small (fixed) rank is 
computable in polynomial time [Barvinok 96]. 

Since the exact computation 
is difficult, the next goal is to find a ``very good'' approximation 
algorithm. A fully polynomial time (randomized) approximation scheme
is a (probabilistic) algorithm that for any given $\epsilon>0$ 
approximates the desired quantity within relative error 
$\epsilon$ in time polynomial in $\epsilon^{-1}$. 
Probabilistic methods based on rapidly mixing Markov chains
resulted in finding such approximation schemes for permanents of 
dense 0-1 matrices (that is, the matrices with at least $n/2$ 1's
in every row and column), random matrices and some special 
0-1 matrices (see [Jerrum and Sinclair 89] and [Jerrum and Sinclair 97]). 
However, for the class of all 0-1 matrices no fully polynomial time 
randomized approximation scheme is known (but there is a 
``mildly exponential'' approximation scheme, see [Jerrum and Vazirani 96]).

In [Barvinok 97b], [Barvinok 99] and [Linial {\it et al.} 20+] a more 
modest goal was posed and achieved. 
Given an arbitrary non-negative $n \times n$ 
matrix $A$, the polynomial time algorithms [Barvinok 97b] and [Barvinok 99] 
(randomized) and [Linial {\it et al.} 20+] (deterministic) produce a 
number $\alpha$ such that 
$$c^n \per A \leq \alpha \leq \per A \tag5.1$$
for some absolute constant $c>0$. Currently the best values of 
$c$ are $c \approx 0.76$ for the randomized algorithm of [Barvinok 99] and 
$c \approx 0.37$ for the deterministic algorithm of [Linial {\it et al.} 20+].
We also note that any polynomial time algorithm achieving a 
subexponential approximation error can be ``upgraded'' to a polynomial 
time approximation scheme, see [Barvinok 99].    

Let $A$ be an $n \times n$ matrix of 0's and 1's. If $\per A$ is 
``big'' (for example, if $\per A$ is of the order $n!/2^n$, which is the 
average value of the permanent for all $n \times n$ 0-1 matrices), 
the additional factor of $c^n$ in (5.1) should not be considered 
as a heavy liability. But if $\per A$ is ``small'' (for example, 
if $\per A$ is of the order $2^{0.01n}$), the lower bound in (5.1) is useless
and the $\alpha$ produced by the algorithms may well be less than 1.
The method developed in this section is designed to provide a partial 
remedy in this situation of a small permanent. 
Our approach should be considered within
the growing family of algorithms that provide a crude 
yet fast and universally applicable estimates. 

Our algorithm for estimating the permanent of a 0-1 matrix $A$ consists 
of constructing a graph $G_A$ as above, finding an economical 
embedding of the set $\FF$ of perfect matchings in $G_A$ into a Boolean cube
(Section 2.4) and estimating the cardinality $|\FF|$ using 
Algorithm 3.3 and Theorem 3.9. We present a summary below.
\subhead (5.2) Algorithm for approximating the permanent \endsubhead

Given an $n \times n$ 0-1 matrix $A=(a_{ij})$, let $G=(V, E)$ be the 
graph with the set of vertices $V=V^+ \cup V^-$, where
$V^+=\{1^+, \ldots, n^+\}$ and $V^-=\{1^-, \ldots, n^-\}$
and the set of edges $E=\bigl\{(i^+, j^-): a_{ij}=1 \bigr\}$.
Let $s_i^+$ be the degree of $i^+$ (the $i$-th row sum of $A$) and 
let $s_j^-$ be the degree of $j^-$ (the $j$-th column sum of $A$).
Let us compute 
$$m^+=\sum_{i=1}^n \lceil \log_2 s_i^+ \rceil \quad  \text{and} 
\quad m^-=\sum_{j=1}^n  \lceil \log_2 s_i^- \rceil$$ 
and let $m=\min(m^+, m^-)$. 

Following Section 2.5, we construct an embedding of the set 
of the perfect matchings in $G$ in $\{0,1\}^m$. 

Without loss of generality we assume that $m=m^+$ (otherwise 
we switch $V^+$ and $V^-$, which corresponds to transposing $A$). Let 
$m_i=\lceil \log_2 s_i^+ \rceil$, so 
$m=m_1 + \ldots +m_n$.

To every edge $e=(i^+, j^-)$ of $G$ incident with $i^+$, let us 
assign a binary string $\phi_i(e)$ of length $m_i$ 
so that $\phi_i(e_1) \ne \phi_i(e_2)$ for every pair of distinct edges 
$e_1$ and $e_2$ with the same endvertex $i^+$. 
 
Given a precision $\epsilon>0$, 
let us generate $k=\lceil 48/m \epsilon^2 \rceil$ 
random binary strings $x_1, \ldots, x_k$, 
each of length $m$.

For each $x=x_i$, $i=1, \ldots, k$, let us do the following procedure:
\medskip
Consider $x$ as a string of $n$ substrings, $x=y_1 \ldots y_n$, where 
$y_i$ is a binary string of length $m_i$. To every edge $e=(i^+, j^-)$ 
of $G$ assign weight $\gamma_e=\dist\bigl(\phi_i(e), y_i\bigr)$, where 
$\dist$ is the Hamming distance between binary strings.
Find the minimum weight $\alpha =\alpha(x)$ of a perfect matching
in $G$ using the
Assignment Problem algorithm, see Section 11.2
of [Papadimitriou and Steiglitz 98].
\medskip    
Compute the average 
$$\alpha={1 \over k} \sum_{i=1}^k \alpha(x_i).$$
Compute 
$$\beta={1 \over 2} -{\alpha \over m}.$$
Output $\beta$.
\proclaim{(5.3) Theorem} Let $A$ be an $n\times n$ 0-1 matrix such 
that $\per A>0$. 

Let $s_1^+, \ldots, s_n^+$ be the row and let 
$s_1^-, \ldots, s_n^-$ be the column sums of $A$. Let 
$$m=\min\biggl\{\sum_{i=1}^n \lceil \log_2 s_i^+ \rceil, \qquad
\sum_{i=1}^n \lceil \log_2 s_i^- \rceil \biggr\}.$$
With probability at least $0.9$, the output $\beta$ of Algorithm 5.2
satisfies 
$$|\beta -\rho| \leq \epsilon,$$
where $0 \leq \rho \leq 1/2$ is a number
such that
$$1-H\Bigl({1 \over 2}-\rho\Bigr) \leq {\log_2 \per A \over m} \leq H(\rho)$$
and $\displaystyle H(x)=x \log_2 {1 \over x} + (1-x) \log_2 {1 \over 1-x}$ 
is the entropy function.
To find $\beta$, Algorithm 5.3 solves $k=\lceil 48/m \epsilon^2 \rceil$ 
Assignment Problems of size $n \times n$.
\endproclaim
\demo{Proof} Let $\FF$ be the set of perfect matchings in the graph
$G=G_A$. The proof follows by the ``economical embedding'' construction 
of Section 2.5, Algorithm 3.3, Theorem 3.6 and Theorem 3.9.
{\hfill \hfill \hfill} \qed
\enddemo
The estimate of Theorem 5.3, however crude, allows us, for example, to decide
in polynomial time whether the permanent of a given $n \times n$ 0-1 
matrix is subexponential in $n$. The precise statement is as follows.
\proclaim{(5.4) Corollary} Let us fix an $0 <\alpha <1$ and let 
us choose any $\beta> (1+\alpha)/2$.
Suppose that $A$ is $n \times n$ 0-1 matrix such that 
$\displaystyle \per A \leq 2^{n^{\alpha}}$. Let us apply Algorithm 5.2 
with $\epsilon=1/m$. Then, for all sufficiently large $n$, the 
estimates of Theorem 5.3 allow us to conclude that  
$\displaystyle \per A \leq 2^{n^{\beta}}$.
\endproclaim
\demo{Proof}   
 We observe that $m \leq n(\log_2 n +1)$. 
By (3.7.2), cf. also Corollary 3.11, we conclude that 
$\rho=O(n^{\alpha/2} m^{-1/2})$ and the proof follows by Theorem 5.3 and 
(3.7.2).
{\hfill \hfill \hfill} \qed
\enddemo
Similarly, one can show that if $\displaystyle \per A \geq 2^{n^{\alpha}}$
then for any $\beta<2\alpha-1$, Algorithm 5.2 with $\epsilon=1/m$ 
would allow us to conclude that $\displaystyle \per A \geq 2^{n^{\beta}}$
for all sufficiently large $n$. The estimate is, of course, void for 
$\beta \leq 1/2$, but it is getting better as $\beta$ approaches 1.
For example, if $\per A$ has the order of $\displaystyle 2^{n^{0.95}}$, 
Algorithm 5.2
would allow us to conclude that $\per A$ is greater than 
$\displaystyle 2^{n^{0.89}}$ and 
is smaller than $\displaystyle 2^{n^{0.98}}$. 

Corollary 5.4 demonstrates something that none of the exponential 
error algorithms (cf. (5.1)) can possibly do (neither can 
any other polynomial time algorithm known to the authors).
On the other hand, algorithms of [Barvinok 97b], [Barvinok 99] and 
[Linial {\it et al.} 20+] are better than Algorithm 5.2 for matrices
with large permanent. Another interesting feature of Algorithm 5.2 is that 
it clearly favors sparse matrices, as the value of $m$
(the dimension of the cubical embedding, see Example 2.5) for such matrices 
is smaller. Algorithms from [Barvinok 97b], [Barvinok 99] and 
[Linial {\it et al.} 20+] seem to be completely indifferent to 
sparseness and even show some inclination to like dense matrices 
better. Thus, in the case of $m=O(n)$ Algorithm 5.2 beats the 
said algorithms on a wider range of permanents (for example, of the 
order $2^{0.01n}$). The final remark is about practical implementation 
of Algorithm 5.2. If $\per A$ is expected to be large enough
(say, of the order $2^{\alpha n}$ for some positive $\alpha$), it suffices
to choose $\epsilon=7m^{-1/2}$, for example. Thus, Algorithm 5.2 boils
down to solving one Assignment Problem. The algorithm should be able 
to handle reasonably sparse matrices with the size $n$ of the order of several
hundreds. 

Our method applies just as well to counting perfect matchings in 
non-bipartite graphs, which is a more general problem. We discussed 
the bipartite case in detail because of its connection with the 
permanent, a problem with rich history and plenty results available 
for comparison.

\head 6. Proofs of Theorems 4.4 and 4.5 \endhead

\definition{(6.1) Definition} We recall that $C_N$ is the Boolean cube 
$\{0,1\}^N$ endowed with the uniform probability measure and that 
$\Lambda_N$ is the Boolean cube $\{0,1\}^N$ endowed with the probability 
measure of Definition 4.1. Let $\Omega_N=C_N \times \Lambda_N$. We 
consider the product measure on $\Omega_N$, so 
$$\Pr\{(x,l)\}=p^{|l|} q^{n-|l|} 2^{-N}, \quad \text{where} \quad 
|l|=\lambda_1 + \ldots + \lambda_N \quad \text{for} \quad l=
(\lambda_1, \ldots,\lambda_N).$$
Hence a point $(x,l) \in \Omega_N$ is interpreted as a vector  
of $2n$ independent random variables 
$(\xi_1, \ldots, \xi_n; \lambda_1, \ldots, \lambda_n)$, where 
$\Pr\{\xi_i=0\}=\Pr\{\xi_i=1\}=1/2$, $\Pr\{\lambda_i=1\}=p$ and 
$\Pr\{\lambda_i=0\}=q$.  
We observe that
$$\Delta(A,p)=\EE d_l(x,A). \leqno(6.1.1)$$ 
\enddefinition
First, we need a version of the concentration inequality (3.4).
\proclaim{(6.2) Lemma} 
Let $A \subset C_N$ be a set.
Then for every $\delta \geq 0$
$$\Pr\Bigl\{(x,l) \in \Omega_N: |d_l(x,A)-\Delta(A,p)| \geq \delta +4 \sqrt{N} 
\Bigr\} \leq 4e^{-\delta^2/N}.$$  
\endproclaim
\demo{Proof}
Given an $A \subset C_N$, let $f: \Omega_N \longrightarrow {\Bbb R}$ be 
defined by $f(x,l)=d_l(x,A)$. Let $M$ be the median of $f$, that is, a
number such that 
$$\Pr\bigl\{(x,l) \in \Omega_N: f(x,l) \leq M \bigr\} \geq 1/2 
\quad \text{and} \quad 
\Pr\bigl\{(x,l) \in \Omega_N: f(x,l) \geq M \bigr\} \geq 1/2.$$  
Since $f$ is a function with 
Lipschitz constant 1, it follows by inequality (2.1.3) of [Talagrand 95] that
$$\Pr\bigl\{(x,l) \in \Omega_N: |f(x,l)-M| \geq 
\delta \bigr\} \leq 4 e^{-\delta^2/N}$$
for any $\delta \geq 0$. 

Since $f$ is integer-valued, we can choose $M$ to be integer. Then
$$\split \EE |f(x,l)-M|&=
\sum_{k=0}^{+\infty} k \Pr\bigl\{(x,l): |f(x,l)-M|=k\bigr\} \\
&=\sum_{k=1}^{+\infty} \Pr\bigl\{(x,l): |f(x,l)-M| \geq k \bigr\} 
\leq 4 \sum_{k=1}^{+\infty} e^{-k^2/N} \\ &\leq 4 \int_0^{+\infty} e^{-x^2/N} 
dx =2\sqrt{\pi N} \leq 4\sqrt{N}.\endsplit$$
Since by (6.1.1) we have
$\Delta(A,p)=\EE f$, we conclude that $|\Delta(A,p)-M| \leq 4 \sqrt{N}$.
Therefore, 
$$\split 
\Pr\bigl\{(x,l): |d_l(x,A)-\Delta(A,p)| 
\geq \delta+ 4\sqrt{N} \bigr\} &\leq 
\Pr\bigl\{(x,l): |d_l(x,A)-M| \geq \delta \bigr\} \\ &\leq
4 \exp\bigl\{-\delta^2/N \bigr\}. \endsplit$$ 
{\hfill \hfill \hfill} \qed
\enddemo
Next, we need an analogue of the scaling trick (3.5).
\proclaim{(6.3) Lemma} Let us fix positive integers $k$ and $n$ and let 
$N=kn$. Let us identify $C_N=(C_n)^k$, $\Lambda_N=(\Lambda_n)^k$ and 
$\Omega_N=(\Omega_n)^k$. Thus a point $(x,l) \in \Omega_N$ is identified 
with $x=(x_1, \ldots, x_k; l_1, \ldots, l_k)$, where $x_i \in C_n$ and 
$l_i \in \Lambda_n$. 

For a subset $A \subset C_n$, let $B=A^k \subset C_N$.
Then 
$$d_l(x, B)=\sum_{i=1}^k d_{l_i}(x_i, A) \quad \text{and} \quad 
\Delta(B,p)=k \Delta(A).$$
\endproclaim
\demo{Proof} Clearly, 
$$d_l(x,y)=\sum_{i=1}^k d_{l_i}(x_i, y_i) \quad \text{for all}
 \quad x,y \in C_N$$
and the first identity follows. Now, by (6.1.1)  
$$\Delta(B,p)=\EE d_l(x, B) =\sum_{i=1}^k \EE d_{l_i}(x_i, A) = k \Delta(A,p).
$$
{\hfill \hfill \hfill} \qed
\enddemo
Now we are ready to prove Theorem 4.4.
\demo{Proof of Theorem 4.4}
 Let $N=nk$ and let us identify 
$C_N=(C_n)^k$, $\Lambda_N =(\Lambda_n)^k$ and 
$\Omega_N=(\Omega_n)^k$.
Let $B=A^k\subset C_N$ as in Lemma 6.3. Applying Lemma 6.2,
we get 
$$\Pr \Bigl\{(x,l) \in \Omega_N: |d_l(x,B)-\Delta(B,p)|
 \geq \delta+4\sqrt{N} \Bigr\} \leq 4 e^{-\delta^2/N}$$
for any $\delta \geq 0$.
Using Lemma 6.3, we conclude:
$$\Pr \Bigl\{(x,l) \in \Omega_N: 
\Big|{1 \over k}\sum_{i=1}^k d_{l_i}(x_i, A)-\Delta(A,p)\Big|
 \geq \delta/k+4\sqrt{n/k} \Bigr\} \leq 4 e^{-\delta^2/N}.$$
Let us choose $\delta=\epsilon k/2$. Hence 
$$\Pr \Bigl\{(x,l) \in \Omega_N: 
\Big|{1 \over k}\sum_{i=1}^k d_{l_i}(x_i, A)-\Delta(A,p)\Big|
 \geq \epsilon/2 +4\sqrt{n/k} \Bigr\} \leq 4 e^{-\epsilon^2 k/n}.$$
Since $k \geq 64n/\epsilon^2$, the proof follows.
{\hfill \hfill \hfill} \qed
\enddemo 
Next, we need a (crude) version of inequality (3.7.1).
\proclaim{(6.4) Lemma} Let $\epsilon \geq 0$, let 
$r(\epsilon)=pN(1-\epsilon)/2$.
Let $y \in C_N$ be a point. Then 
$$\Pr\bigl\{(x,l) \in \Omega_N: d_l(x,y) \leq r(\epsilon) \bigr\} 
\leq e^{-\epsilon^2 p N/4}.$$
\endproclaim
\demo{Proof} Without loss of generality we may assume that $y=0$. 
Then 
$$\Pr\bigl\{(x,l) \in \Omega_N: d_l(x,0) \leq r(\epsilon) \bigr\}=
\Pr\Bigl\{(x,l) \in \Omega_N: \sum_{i=1}^N \xi_i \lambda_i \leq 
r(\epsilon) \Bigr\},$$
where $x=(\xi_1, \ldots, \xi_N)$ and $l=(\lambda_1, \ldots, \lambda_N)$.
Let $\zeta_i=\xi_i \lambda_i$. Then $\zeta_i$, $i=1, \ldots, N$ are 
independent random variables such that 
$\Pr\{\zeta_i=1\}=p/2$ and $\Pr\{\zeta_i=0\}=1-p/2$.
Hence
$$\Pr\bigl\{(x,l) \in \Omega_N: d_l(x,y) \leq r(\epsilon)\bigr\} = 
\Pr\bigl\{\zeta_1 + \ldots  + \zeta_N \leq r(\epsilon) \bigr\} 
\leq e^{-\epsilon^2 p N/4}$$
by a corollary of Chernoff's inequality (see [McDiarmid 89]).
{\hfill \hfill \hfill} \qed
\enddemo
Now we are ready to prove the first part of Theorem 4.5.
\demo{Proof of inequality (4.5.1)}
Let us choose a positive integer $m$, let $N=mn$, 
let $C_N=(C_n)^m$, and let $\Lambda_N=(\Lambda_n)^m$. Let 
$B=A^m\subset C_N$ as in Lemma 6.3. 

Let us choose an $\alpha>0$. Applying Lemma 6.4, we obtain
$$\Pr\Bigl\{(x,l) \in \Omega_N: d_l(x,B) \leq pN(1-\sqrt{\alpha})/2 \Bigr\}
 \leq 
|B| e^{-\alpha p N/4}=(|A| e^{-\alpha p n/4})^m.$$
Therefore, by Lemma 6.3
$$\Pr\Bigl\{(x,l) \in \Omega_N: {1 \over m} \sum_{i=1}^m d_{l_i} 
(x_i, A) \leq pn(1-\sqrt{\alpha})/2 \Bigr\} \leq (|A| e^{-\alpha p n/4})^m.$$
The right hand side of the inequality tends to $0$ provided
$\alpha > 4 \ln |A|/pn$.
Since by the Law of Large Numbers
$${1 \over m} \sum_{i=1}^m d_{l_i} (x_i, A) \longrightarrow \Delta(A,p)
\quad \text{in probability} \quad \text{as} \quad
m \longrightarrow +\infty,$$ we must have 
$$\Delta(A,p) \geq pn(1-\sqrt{\alpha})/2 \quad \text{for any} \quad
\alpha > 4 \ln |A|/pn.$$
Hence 
$$\Delta(A,p) \geq pn(1-\sqrt{\alpha})/2 \quad \text{for} \quad
\alpha = 4 \ln |A|/pn,$$
which is equivalent to (4.5.1).
{\hfill \hfill \hfill} \qed
\enddemo
In Section 3, we used the sharp isoperimetric inequality (Theorem 3.8) for 
the Hamming distance in $C_n$ to get a sharp upper bound for 
$n^{-1} \log_2 |A|$.
Unfortunately, we don't know of a similar result for the randomized 
Hamming distance. To prove (4.5.2)--(4.5.3), we proceed by induction on $n$ 
in a way resembling that of [Talagrand 95] (see also Remark 6.9).  

We start with a simple technical result.
\proclaim{(6.5) Lemma} For any $0 \leq \epsilon \leq 1$, any 
$\gamma \geq 0$ and any $0 < p \leq 1$ and $q=1-p$ we have 
$$\min\Bigl\{{p \gamma \over 2} +\ln {1 \over 1+\epsilon}, \quad 
p \ln {1 \over 1-\epsilon} + q \ln {1 \over 1+\epsilon} \Bigr\} \leq 
\max\Bigl\{0,\quad \ln(1+e^{\gamma/2})-{q\gamma \over 2} - \ln 2 \Bigr\}.$$
\endproclaim
\demo{Proof} Fixing $p,q$ and $\gamma$, let 
$$f(\epsilon)={p \gamma \over 2} +\ln {1 \over 1+\epsilon} \quad \text{and} 
\quad g(\epsilon)=p \ln {1 \over 1-\epsilon} + q \ln {1 \over 1+\epsilon}.$$
Then $f(0) \geq 0$ and $f(\epsilon)$ is decreasing
whereas $g(\epsilon)$ behaves as follows: $g(0)=0$ and if $p \geq q$ then 
$g(\epsilon)$ is increasing and if $p<q$ then $g(\epsilon)$ is 
decreasing for $0<\epsilon<q-p$ and increasing for $q-p < \epsilon <1$.
Furthermore, 
$f(\epsilon_0)=g(\epsilon_0)$ at the single point 
$\epsilon_0=(e^{\gamma/2}-1)/(1+e^{\gamma/2})$, where
$f(\epsilon_0)=g(\epsilon_0)=\ln(1+e^{\gamma/2})-q\gamma/2 - \ln 2$.
The proof now follows.
{\hfill \hfill \hfill} \qed
\enddemo  
\definition{(6.6) Definition} Let 
$\mu_n$ (or simply $\mu$) denote the uniform probability measure in $C_n$.
Hence $\mu(A)=|A|/2^n$.
\enddefinition
The induction is based on the following lemma.
\proclaim{(6.7) Lemma} Let $A \subset C_{n+1}$ be a set.
Let 
$$A_0=\bigl\{x \in C_n: (x,0) \in A \bigr\} \quad \text{and} 
\quad A_1=\bigr\{x \in C_n: (x,1) \in A \bigr\}.$$
For $l \in \Lambda_n$ let $(l,0) \in \Lambda_{n+1}$ denote 
$l$ appended by $\lambda_{n+1}=0$ and let $(l,1) \in \Lambda_{n+1}$ denote
$l$ appended by $\lambda_{n+1}=1$. Let  
$$\Delta_0(A,p)=\EE d_{(l,0)}(x, A) \quad \text{and} \quad 
\Delta_1(A,p)=\EE d_{(l,1)}(x, A),$$
where the expectation is taken with respect to a random 
$(x,l) \in C_{n+1} \times \Lambda_n$. 
Then 
$${\mu_n(A_0) + \mu_n(A_1) \over 2} =\mu_{n+1}(A); \leqno(6.7.1)$$
$$\Delta(A,p)=q\Delta_0(A,p) + p\Delta_1(A,p); \leqno(6.7.2)$$
$$\Delta_0(A,p) \leq \Delta(A_i,p) \quad \text{for} \quad i=0,1; 
\leqno(6.7.3)$$
$$\Delta_1(A,p) \leq \Delta(A_i,p)+{1 \over 2} \quad \text{for} \quad 
i=0,1; \leqno(6.7.4)$$
$$\Delta_1(A,p) \leq {\Delta(A_0,p) + \Delta(A_1,p) \over 2}. \leqno(6.7.5)$$
\endproclaim
\demo{Proof} Clearly, $|A_0|+|A_1|=|A|$, so (6.7.1) follows.
Identity (6.7.2) is immediate from Definitions 4.1.
We observe that for any $x, y \in C_n$,
$$d_{(l,0)}((x,j), (y,i))=d_l(x,y), \quad \text{where} \quad i,j \in \{0,1\}.$$
Hence
$$d_{(l,0)}((x,j), A) \leq d_l(x, A_i), \quad i,j=0,1$$ 
and (6.7.3) follows by averaging.

Next, we observe that 
$$d_{(l,1)}((x,i), (y,j))=\cases d_l(x,y) &\text{if\ }i=j \\ d_l(x,y)+1
&\text{if\ } i \ne j. \endcases$$
Therefore, 
$$d_{(l,1)}((x,1), A)=\min\bigl\{d_l(x, A_1), 
d_l(x, A_0)+1 \bigr\}$$
and
$$d_{(l,1)}((x,0), A)=\min\bigl\{d_l(x, A_0), 
d_l(x, A_1)+1 \bigr\}.$$
Averaging over $(x,l) \in C_{n+1} \times \Lambda_n$, we get
$$\split \Delta_1(A,p)
 &=\EE d_{(l,1)}(x,A)={\EE d_{(l,1)}\bigl((x,1), A\bigr) + 
\EE d_{(l,1)}\bigl((x,0), A\bigr) \over 2} \\ &\leq 
{\EE d_l(x, A_1)+ \EE d_l(x, A_1) +1 \over 2} = \Delta(A_1, p) + {1 \over 2}.
\endsplit$$
Similarly,
$$\split \Delta_1(A,p)
 &=\EE d_{(l,1)}(x,A)={\EE d_{(l,1)}\bigl((x,1), A\bigr) + 
\EE d_{(l,1)}\bigl((x,0), A\bigr) \over 2} \\ &\leq 
{\EE d_l(x, A_0)+ 1+ \EE d_l(x, A_0) \over 2} = \Delta(A_0, p) + {1 \over 2},
\endsplit$$
which completes the proof of (6.7.4).
Finally, 
$$\split \Delta_1(A,p)
 &=\EE d_{(l,1)}(x,A)={\EE d_{(l,1)}\bigl((x,1), A\bigr) + 
\EE d_{(l,1)}\bigl((x,0), A\bigr) \over 2} \\ &\leq 
{\EE d_l(x, A_1)+ \EE d_l(x, A_0)\over 2} = {\Delta(A_1, p) +\Delta(A_0,p) 
\over 2}\endsplit$$
and (6.7.5) is proved.
{\hfill \hfill \hfill} \qed
\enddemo
Now we use induction to get a preliminary bound.
\proclaim{(6.8) Lemma} Suppose that for some $\gamma \geq 0$, 
$0< p \leq 1$ and $q=1-p$,
$$\ln \bigl(1+e^{\gamma/2}\bigr) - {q \gamma \over 2} -\ln 2 \geq 0.$$
Then for any non-empty set $A \subset C_n$ we have 
$$\gamma \Delta(A,p) +\ln \mu(A) \leq n\Bigl( \ln \bigl(1+e^{\gamma/2} \bigr)-
{q \gamma \over 2} -\ln 2\Bigr).$$
\endproclaim
\demo{Proof} We proceed by induction on $n$. 
If $n=1$ then the two cases are possible:
\medskip
$A$ consists of a single point, $\mu(A)=1/2$ and $\Delta(A,p)=p/2$ (see 
Example 4.2);
\smallskip
$A=\{0,1\}$, $\mu(A)=1$ and $\Delta(A,p)=0$.
\medskip
In both cases the inequality holds.

Suppose that the inequality holds for non-empty subsets of 
$C_n$. Let us prove
that it holds for non-empty $A \subset C_{n+1}$.  Let us define 
$A_0, A_1 \subset C_n$ as in Lemma 6.7. From (6.7.1)
it follows that either 
$$\mu_n(A_0)=(1-\epsilon) \mu_{n+1}(A) \quad \text{and} \quad 
\mu_n(A_1)=(1+\epsilon) \mu_{n+1}(A)$$ or 
$$\mu_n(A_1)=(1-\epsilon) \mu_{n+1}(A) \quad \text{and} \quad  
\mu_n(A_0)=(1+\epsilon) \mu_{n+1}(A)$$ for some $0 \leq \epsilon \leq 1$.

Let $B$ be the one of the sets $A_0$, $A_1$ that has a bigger measure 
$\mu_n$ (either of the two if $\mu_n(A_0)=\mu_n(A_1)$) and let
$D$ be the one of the sets $A_0$, $A_1$ that has a bigger value of 
$\Delta(\cdot, p)$ (either of the two if $\Delta(A_0,p)=\Delta(A_1, p)$).
Then 
$$\mu_n(B) \geq (1+\epsilon) \mu_{n+1}(A) \quad \text{and} \quad
\mu_n(D) \geq (1-\epsilon) \mu_{n+1}(A).$$
Furthermore, by (6.7.3)
$$\Delta_0(A,p) \leq \Delta(B,p) \quad \text{and} \quad
\Delta_0(A,p) \leq \Delta(D, p)$$
whereas by (6.7.3) and (6.7.5) 
$$\Delta_1(A, p) \leq \Delta(B,p) +{1 \over 2} \quad \text{and} 
\quad \Delta_1(A,p) \leq \Delta(D,p).$$  
Hence we get
$$\gamma \Delta_0(A,p) +\ln \mu_{n+1}(A) \leq 
\gamma \Delta(B,p) +\ln \mu_n(B) + \ln {1 \over 1+\epsilon}$$ 
and 
$$\split & \gamma \Delta_1(A,p) +\ln \mu_{n+1}(A) \leq 
\\ 
&\min\Bigl\{ \gamma \Delta(B,p) + \ln \mu_n(B) + \ln {1 \over 1+\epsilon}
+{\gamma \over 2}, \quad \gamma \Delta(D,p) +\ln \mu_n(D) + 
\ln {1 \over 1-\epsilon} \Bigr\}. \endsplit $$
Clearly, $B$ is non-empty. Assume first, that $D$ is non-empty as well. 
Applying the induction hypothesis to $B$ and $D$, we conclude
that
$$\gamma \Delta_0(A,p) +\ln \mu_{n+1}(A) \leq 
 n\Bigl( \ln \bigl(1+e^{\gamma/2} \bigr)-
{q \gamma \over 2} -\ln 2\Bigr) + \ln {1 \over 1+\epsilon} $$ and 
$$\gamma \Delta_1(A,p) +\ln \mu_{n+1}(A) \leq
 n\Bigl( \ln \bigl(1+e^{\gamma/2} \bigr)-
{q \gamma \over 2} -\ln 2\Bigr) + \min\Bigl\{\ln {1 \over 1+\epsilon}
+{\gamma \over 2}, \quad  
\ln {1 \over 1-\epsilon} \Bigr\}.$$
Adding the first inequality multiplied by $q$ and the second inequality 
multiplied by $p$ and using (6.7.2), we get 
$$\split &\gamma \Delta(A,p) + \ln \mu_{n+1}(A)  \leq \\
&n\Bigl( \ln \bigl(1+e^{\gamma/2} \bigr)-
{q \gamma \over 2} -\ln 2\Bigr)+ 
\min\Bigl\{{p \gamma \over 2} +\ln {1 \over 1+\epsilon}, \quad
p \ln {1 \over 1-\epsilon} + q \ln {1 \over 1+\epsilon} \Bigr\}.
\endsplit$$
The desired inequality follows by Lemma 6.4.

If $D$ is empty then $\mu_n(B)=2\mu_{n+1}(A)$ and we obtain
$$\gamma \Delta_0(A,p) + \ln \mu_{n+1}(A) \leq \gamma \Delta(B,p) + 
\ln \mu_n(B) - \ln 2$$ and 
$$\gamma \Delta_1(A,p) + \ln \mu_{n+1}(A) \leq \gamma \Delta(B,p) + 
\ln \mu_n(B) - \ln 2 + {\gamma \over 2}$$
Adding the first inequality multiplied by $q$ to the second inequality 
multiplied by $p$ and using (6.7.2) and the induction hypothesis, we 
get:
$$\split \gamma \Delta(A,p)+ \ln \mu_{n+1}(A) &\leq \gamma \Delta(B,p) +
\ln \mu_n(B) - \ln 2 +{p\gamma \over 2} \\
 &\leq n\Bigl( \ln \bigl(1+e^{\gamma/2} \bigr)-
{q \gamma \over 2} -\ln 2\Bigr) + 
\Bigl({\gamma \over 2} - {q\gamma \over 2} -\ln 2\Bigr) \\ &\leq 
 (n+1)\Bigl( \ln \bigl(1+e^{\gamma/2} \bigr)-
{q \gamma \over 2} -\ln 2\Bigr), \endsplit$$ 
which completes the proof.
{\hfill \hfill \hfill} \qed
\enddemo
Now we are ready to complete the proof of Theorem 4.5.
\demo{Proof of (4.5.2)--(4.5.3)} 
By Lemma 6.8, 
$${\ln |A| \over n}={\ln \mu_n(A) \over n} +\ln 2 \leq 
\ln(1+e^{\gamma/2})-{q\gamma \over 2} - {\gamma \Delta(A,p) \over n}=
\ln(1+e^{\gamma/2})-{\gamma \over 2} +\gamma \rho$$
provided 
$$\ln \bigl(1+e^{\gamma/2}\bigr) - {q \gamma \over 2} -\ln 2 \geq 0.$$
We optimize the inequality on $\gamma \geq 0$.
Let
$$\gamma=2 \ln \Bigl({1 \over 2\rho}-1 \Bigr).$$
Since we assumed that $\rho \leq 1/4$, we have $\gamma \geq 0$.
Furthermore,
$$\split \ln \bigl(1+e^{\gamma/2}\bigr) - {q \gamma \over 2} -\ln 2 &=
\ln {1 \over 2 \rho} - q \ln \Bigl({1 \over 2\rho}-1 \Bigr) - \ln 2 \\
&=-\ln(1-2\rho) +p \bigl( \ln(1-2\rho) -\ln (2\rho) \bigr) -\ln 2 \geq 0,
\endsplit$$ 
because of (4.5.2).
Therefore,
$${\ln |A| \over n} \leq \ln {1 \over 2\rho} - \ln {1-2\rho \over 2\rho} 
+ 2\rho \ln {1 -2\rho \over 2\rho}=2 \rho \ln {1 \over 2 \rho} +
(1-2\rho) \ln {1 \over 1-2\rho}$$
and (4.5.3) follows.
{\hfill \hfill \hfill} \qed
\enddemo
\subhead (6.9) Remark \endsubhead
Our proof of (4.5.2)--(4.5.3) can be considered 
as an ``additive'' version of Talagrand's method [Talagrand 95]. Indeed,
Talagrand's approach very roughly can be can be stated as follows.
Let $\Omega$ be a space with the distance function $d$ and 
probability measure $\mu$. To prove an isoperimetric inequality 
for $A \subset \Omega$, we first find a uniform bound for
the expression $\mu^{\alpha}(A) \cdot \EE \exp\{\tau d(x,A)\}$ and then 
adjust parameters $\alpha>0$ and $\tau>0$. This way 
tight inequalities are obtained in [Talagrand 95] for sets $A$ of large 
measure, most often with $\mu(A) \geq 1/2$.
We are mostly interested in sets of a small measure. One can check that 
for ``small sets'' $A$ the inequalities of [Talagrand 95] 
are very far from sharp,
which is, of course, should not be perceived as a ``fault'' of the 
method, since the method was designed for totally different problems. 
We find a uniform bound for the expression $\ln \mu(A) + \alpha 
\EE d(x,A)$, which looks like Talagrand' functional with ``exp'' removed.
Our method seems to produce reasonably good bounds for small sets 
$A$ but it fails miserably for large $A$, with $\mu(A) =1/2$, say. 
As should have been expected, the case of ``middle-sized'' sets is the most
complicated.  

\head 7. Concluding Remarks \endhead

\remark{Connections to Monte-Carlo methods} The main idea of our 
approach can be described as follows: given a (finite) ambient space 
$\Omega$ and a set $A \subset \Omega$, we estimate the cardinality $|A|$ 
by choosing a certain distance function $d$ in $\Omega$ and estimating 
the average distance 
$$\Delta(A)={1 \over |\Omega|} \sum_{x \in \Omega} d(x,A), \quad 
\text{where} \quad d(x,A)=\min_{y \in A} d(x,y)$$
from $x \in \Omega$ to $A$. 
We get the classical Monte-Carlo method if the distance function $d$ 
is the simplest possible:
$$d(x,y)=\cases 1 &\text{if\ } x\ne y \\ 0 &\text{if \ } x=y. \endcases$$
In this case, $\Delta(A)=|A|/|\Omega|$, so there is a direct relation 
between $\Delta(A)$ and $|A|$. It is well understood that the main 
difficulty with the Monte-Carlo method is that if $|A|$ is very small 
compared to $|\Omega|$, it is hard to get an estimate for the cardinality 
of $A$ different from 0. In other words, if $|A|$ is 
``exponentially small'' compared to $|\Omega|$, to get a non-trivial 
bound for $|A|$, we have to compute $\Delta(A)$ with exponentially high 
precision. In this paper, we showed that in many interesting cases 
one can choose a different distance function $d$, so that 
the distance $d(x, A)$ from a point $x \in \Omega$ to $A$ is efficiently 
computable and to get a meaningful estimate of $|A|$ even for 
exponentially small sets $A$, one need to compute $\Delta(A)$ with 
a polynomial precision.

Hence our approach can be considered as a natural extension of 
the Monte-Carlo method. In this context, economical embedding 2.4 can be 
considered as an analogue of the ``importance sampling'', whose 
objective is to replace a large ambient space $\Omega$ by a smaller 
space containing $A$. 
\endremark
\remark{Embedding in different metric spaces} Given a combinatorially 
defined family $\FF \subset 2^X$, we constructed its embedding into the
Boolean cube $\{0,1\}^n$ and investigated what happens in the cube 
is endowed either with the standard Hamming distance (Section 3) or 
with its randomized version (Section 4). In many cases, there are 
different ways of metrization of $\FF$. One example is provided by the 
set $\FF$ of perfect matchings in a given bipartite graph studied 
in the paper.

Let $G=(V^+ \cup V^-, E)$ be a bipartite graph with 
$V^+=\{1^+, \ldots, n^+\}$, $V^-=\{1^-, \ldots, n^-\}$ (cf. Example 1.2
and Section 5). For every vertex $i^+ \in V_+$, let 
$$\Omega_i=\bigl\{e=(i^+, j^-): e \in E \bigr\}$$
be the set of edges of $G$ coming out of $i^+$. Let 
$${\bold \Omega} =\Omega_1 \times \ldots \times \Omega_n.$$
Every perfect matching in $G$ can be identified with a point in 
${\bold \Omega}$, 
so the set $\FF$ of all perfect matchings in $G$ is identified with 
a subset $F \subset {\bold \Omega}$.

Let $d_i$ be a distance function on $\Omega_i$, $i=1, \ldots, n$.
Let us define the distance function $d$ on ${\bold \Omega}$ by 
$$d(x,y)=\sum_{i=1}^n d_i(x_i, y_i), \quad \text{where} \quad 
x=(x_1, \ldots, x_n) \quad \text{and} \quad y=(y_1, \ldots, y_n).$$
It is easy to check that for any $x \in {\bold \Omega}$,
$x=(x_1, \ldots, x_n)$ the distance $d(x, F)$ is the minimum weight 
of a perfect matching in $G$ with weighting $\gamma(e)=d_i(e, x_i)$ 
for $e=(i^+, j^-)$. Hence for any $x$, the value of $d(x, F)$ can be found 
in $O(n^3)$ time.  How should we choose $d_i$ to get the best possible
estimates for the number $|F|$ of perfect matchings in $G$? 

The authors looked into some of the most obvious candidates, when $d_i$ is a 
graph metric on $\Omega_i$ for a complete graph and for a path (circle).
Interestingly, choosing $\Omega_i$ isometric to a 
subset of a power of the complete graph with 
$m=\Theta(\ln n)$ vertices leads to an improvement by 
logarithmic factor similar to that of the randomized Hamming 
distance (Section 4). The choice of $d_i$ used in this 
paper comes from identifying $\Omega_i$ with a subset of the Boolean 
cube $\{0,1\}^{m_i}$ for $m_i=\lceil \log_2 |\Omega_i| \rceil$, see 
Section 5. Perhaps one should use a whole family of distance functions 
$d_i$ and combine the resulting estimates. General isoperimetric 
inequality of [Alon {\it et al.} 98] may be very useful for that.
\endremark
\remark{Weighted counting} Let $\FF \subset 2^X$ be a family of subsets 
of the ground set $X=\{1, \ldots, n\}$ and let 
$\mu(i)=p_i/q_i>0$ be a rational weight of $i \in X$, where
$p_i, q_i \in {\Bbb N}$. 
Let us define
$$\mu(Y)=\prod_{i \in Y} \mu(i) \quad \text{for} \quad 
Y \in \FF \qquad \text{and} \qquad \mu(\FF)=\sum_{Y \in \FF} \mu(Y).$$  
We may be interested to estimate $\mu(\FF)$.
There are several ways to extend our methods to problems of this type, 
here we sketch one.
For every $i \in X$, let 
$m_i=\lceil \log_2 (p_i +q_i) \rceil $. 
Let us choose subsets $A_i \subset C_{m_i}$ and $B_i \subset  C_{m_i}$ 
such that $|A_i|=p_i$, $|B_i|=q_i$ and $A_i \cap B_i=\emptyset$. 
Let $m=m_1 + \ldots + m_n$ and let us 
identify 
$$C_m=C_{m_1} \times \ldots \times C_{m_n}.$$
For $Y \subset \FF$ let $Z_Y \subset C_m$ be the direct product 
of $n$ factors, the $i$-th factor being $A_i$ if $i \in Y$ and $B_i$ if 
$i \notin Y$. Finally, let $F \subset C_m$ be the union of all $Z_Y$ 
for $Y \in \FF$. We see that $\mu(\FF)=(q_1 \cdots q_n)^{-1} |F|$. 
Moreover, one can 
define subsets $A_i$ and $B_i$ in such a way that Optimization Oracle 1.1 for 
$\FF$ gives rise to Distance Oracle 2.2 for $F$. 
This construction corresponds to the straightforward embedding (2.3). In some 
cases, there is a way to come up with an economical embedding in the spirit 
of (2.4). 
\endremark

\head Acknowledgment \endhead

The authors are grateful to M. Gromov and B. Sudakov 
for many helpful discussions.

\head References \endhead

\noindent [Alon {\it et al.} 98] N. Alon, R. Boppana and J. Spencer,
An asymptotic isoperimetric 

inequality, {\it Geom. Funct. Anal.} 
{\bf 8} (1998), 411--436. 
\smallskip
\noindent [Barvinok 96] A. Barvinok, Two algorithmic results for the traveling 
salesman 

problem, {\it Math. Oper. Res.} {\bf 21} (1996), 65--84.
\smallskip
\noindent [Barvinok 97a] 
A. Barvinok, Approximate counting via random optimization,

{\it  Random Structures $\&$ Algorithms}, {\bf 11}(1997), no. 2, 187--198.
\smallskip
\noindent [Barvinok 97b] A. Barvinok, 
Computing mixed discriminants, mixed volumes, and 

permanents,
{\it Discrete Comput. Geom.} {\bf 18} (1997), no. 2, 205--237.
\smallskip
\noindent [Barvinok 99] 
A. Barvinok, Polynomial time algorithms to approximate permanents 

and mixed
discriminants within a simply exponential factor, 

{\it Random Structures  $\&$ Algorithms}, {\bf 14}(1999), no. 1, 29--61. 
\smallskip
\noindent [Gallucio and Loebl 99] A. Galluccio and M. Loebl,
On the theory of Pfaffian 

orientations. I. Perfect matchings and permanents,
{\it Electron. J. Combin.} {\bf 6}

(1999), no. 1, Research Paper 6, 18 pp. (electronic).
\smallskip
\noindent[Grigoriev and Karpinski 87], D. Yu Grigoriev and M. Karpinski, 
The matching 

problem for bipartite graphs with polynomially bounded 
permanents is in NC, 

{\it Proc. Twenty-Eight Annual IEEE Symp. 
Foundations of Computer Science}, 

IEEE, Computer Society Press, 
Washington, DC, 1987, 162--172.
\smallskip
\noindent [Jerrum 95] M. Jerrum, The computational complexity of counting, 

{\it  Proceedings of the International Congress of Mathematicians}, Vol. 1,2 

(Z\"urich, 1994), 1407--1416, Birkh\"auser, Basel, 1995.
\smallskip
\noindent [Jerrum and Sinclair 89] 
M. Jerrum and A. Sinclair, Approximating the permanent, 

{\it SIAM J. Comput.}, {\bf 18} (1989), no. 6,
1149--1178.
\smallskip
\noindent [Jerrum and Sinclair 97] M. Jerrum and A. Sinclair, The 
Markov chain Monte Carlo

Method: an approach to approximate counting 
and integration, in: 

{\it Approximation Algorithms for NP-hard Problems},
ed. D.S. Hochbaum, PWS, 

Boston, 1997, 483--520. 
\smallskip
\noindent [Jerrum {\it et al.} 86] 
M. Jerrum, L.G. Valiant and V.V. Vazirani, Random generation 

of combinatorial
structures from a uniform distribution, {\it Theoret. Comput. Sci.},

{\bf 43} (1986), no. 2-3, 169--188.
\smallskip
\noindent [Jerrum and Vazirani 96] M. Jerrum and U. Vazirani, 
A mildly exponential 

approximation algorithm for the permanent,
{\it Algorithmica}, {\bf 16}(1996), no. 4-5, 

392--401.
\smallskip
\noindent [Leader 91] I. Leader, Discrete isoperimetric inequalities, in: 
{\it  Probabilistic 

Combinatorics and its Applications 
(San Francisco, CA, 1991)}, 57--80, Proc. 

Sympos. Appl. Math., 44, 
Amer. Math. Soc., Providence, RI, 1991.
\smallskip
\noindent [Linial {\it et al.} 20+] 
N. Linial, A. Samorodnitsky and A. Wigderson, A deterministic 

strongly 
polynomial algorithm for matrix scaling and approximate permanents,

{\it Combinatorica}, to appear.
\smallskip
\noindent [van Lint 99] J.H. van Lint, {\it Introduction to Coding Theory},
Third edition. 

Graduate Texts in Mathematics, {\bf 86}, Springer-Verlag, Berlin, 1999.
\smallskip
\noindent [Lov\'asz and Plummer 86]
L. Lov\'asz and M.D. Plummer, {\it Matching Theory},

North-Holland Mathematics Studies, {\bf 121}, Annals of
Discrete Mathematics, {\bf 29}, 

North-Holland Publishing Co., 
Amsterdam-New York; Akad\'emiai Kiad\'o 

(Publishing House of
the Hungarian Academy of Sciences), Budapest, 1986.
\smallskip
\noindent[McDiarmid 89] C. McDiarmid, On the method of bounded differences,
in: 

{\it Surveys in Combinatorics}, 1989 (Norwich, 1989), 148--188, 
London Math.
Soc.

 Lecture Note Ser., {\bf 141}, Cambridge Univ. Press, Cambridge, 1989
\smallskip
\noindent [Milman and Schechtman 86]
V.D. Milman and G. Schechtman, 
{\it Asymptotic Theory} 

{\it  of Finite-Dimensional Normed Spaces.
With an Appendix by M. Gromov}, 

Lecture Notes in Mathematics, {\bf 1200}, Springer-Verlag, Berlin-New York, 
1986. 
\smallskip
\noindent [Minc 78] H. Minc, {\it Permanents}, 
Encyclopedia of Mathematics and its Applications, 

{\bf 6},
Addison-Wesley, Reading, Mass., 1978.
\smallskip
\noindent [Papadimitriou and Steiglitz 98] 
C.H. Papadimitriou and K. Steiglitz, 

{\it Combinatorial Optimization:
Algorithms and Complexity}, Dover, NY, 1998.
\smallskip
\noindent [Talagrand 95] 
M. Talagrand, Concentration of measure and isoperimetric 

inequalities in 
product spaces, {\it Inst.
Hautes \'Etudes Sci. Publ. Math.} No. 81, 

(1995), 73--205.
\smallskip
\noindent [Valiant 79] L.G. Valiant, The complexity of computing the permanent,

{\it Theoret. Comput. Sci.}, {\bf  8 }(1979), no. 2,
189--201.
\enddocument

\end